\documentclass{elsart}
\usepackage{graphicx,psfrag}
\usepackage{amsmath}
\usepackage{amssymb}
\begin{document}
\begin{frontmatter}

\title{Geometric integrators for 
 multiplicative viscoplasticity: analysis of error accumulation }
\author{A.V. Shutov\corauthref{cor}},
\corauth[cor]{Corresponding author. Tel.: +49-0-371-531-35024; fax: +49-0-371-531-23419.}
\ead{alexey.shutov@mb.tu-chemnitz.de}
\author{R. Krei{\ss}ig}

\address{Institute of Mechanics,
Chemnitz University of Technology,
Str. d. Nationen 62, D-09111 Chemnitz, Germany}

\begin{abstract}
The inelastic incompressibility is a typical feature of metal
plasticity/viscoplasticity. Over the last decade, there has
been a great amount of research related to construction of
numerical integration algorithms which
exactly preserve this geometric property. In this paper we examine,
both numerically and mathematically, the excellent accuracy and
convergence characteristics of such geometric integrators.

In terms of a classical model of finite viscoplasticity, we illustrate
the notion of exponential stability of the exact solution. We show that
this property enables the construction of
effective and stable numerical algorithms, if incompressibility
is exactly satisfied. On the other hand, if the incompressibility
constraint is violated,
spurious degrees of freedom are introduced. This results in the loss
of the exponential stability and a dramatic deterioration of convergence behavior.

\end{abstract}
\end{frontmatter}


\section*{Nomenclature}

\begin{tabbing}
$\mathbf C_{\text{i}}$ \quad \quad \quad \quad \quad  \quad \quad  \= inelastic right Cauchy-Green tensor (see \eqref{intIntv}) \\
$\tilde{\mathbf T}$ \> 2nd Piola-Kirchhoff tensor (see \eqref{prob2}) \\
$\mathbf 1$ \> second-rank identity tensor \\
$\mathbf A \cdot \mathbf B = \mathbf A \mathbf B$ \> product (composition) of two second-rank tensors \\
$\mathbf A : \mathbf B$ \> scalar product of two second-rank tensors \\
$\mathbf A \otimes \mathbf B$ \> tensor product of two second-rank tensors \\
$\| \mathbf A \|$ \> $l_2$ norm of a second-rank tensor (Frobenius norm) \\
$\| \mathbf A \|^*$ \> induced norm of a second-rank tensor (spectral norm)  (see \eqref{Openo}) \\
$(\cdot)^{\text{D}}$ \> deviatoric part of a tensor \\
$ (\cdot)^{\text{T}}$  \> transposition of a tensor \\
$ (\cdot)^{-\text{T}}$ \>  inverse of transposed \\
$\text{tr}(\cdot)$ \> trace of a second-rank tensor \\
$\overline{(\cdot)}$ \> unimodular part of a tensor (see \eqref{UniPro}) \\
$\text{sym}(\cdot)$ \> symmetric part of a tensor   \\
$\langle x \rangle$ \> MacCauley bracket (see $(\ref{examp3})_3$) \\
$\psi_{\text{el}}$ \>  specific free energy density \\
$\text{Dist}(\cdot, \cdot)$ \>  "distance" between two solutions (see \eqref{EnerDist1}) \\
$K$ \> yield stress \\
$\lambda_{\text{i}}$ \> proportionality factor (inelastic multiplier)  (see $(\ref{prob3})_1$) \\
$f$ \> overstress  (see $(\ref{prob3})_2$) \\
$\mathfrak{F}$ \> norm of the driving force (see $(\ref{prob3})_3$) \\
$Sym$ \> space of symmetric second-rank tensors \\
$M$ \> invariant manifold (cf. \eqref{manifold}, \eqref{geopr0}) \\
$\rho_{\scriptscriptstyle \text{R}}$ \> mass density in the reference configuration \\
$k$ \> bulk modulus (see \eqref{spec1})  \\
$\mu$ \>  shear modulus (see \eqref{spec1}) \\
\end{tabbing}

\section{Introduction}

The mechanical processing of materials may involve very large inelastic deformations.
For instance, for equal channel angular extrusion of aluminum alloys,
the introduced accumulated inelastic strain usually varies between 100 and 900 Percent
(depending on the number of extrusions \cite{Valiev}). Even larger deformations can
be introduced by some incremental forming procedures like
spin extrusion \cite{Michel} (the accumulated inelastic strain ranges up to 1000 Percent).
Due to the highly nonlinear character of the underlying mechanical problem,
a correct numerical simulation of such "long" processes is by no means a trivial task.
It is desirable to have numerical algorithms which would be stable with respect
to numerical errors, even if working with big time intervals and big time steps.

The assumption of
\emph{exact inelastic incompressibility} is widely implemented for construction
of material models of metal plasticity and creep (see, for instance, \cite{Haupt}).
Extensive studies were carried out concerning
the construction of numerical integration algorithms which
exactly preserve the incompressibility of the inelastic flow
\cite{DettRes, Hartmann2, Helm2, Luhrs, Miehe, Shutov1, Simo, Svendsen}.\footnote{
The incompressibility condition is given by a linear
invariant in the case of infinitesimal strains inelasticity. Since the linear invariants
are exactly conserved by most of integration procedures (cf. \cite{Hair2}), the problem of
the conservation of incompressibility only appears when working with finite strains.}

In this paper, we asses those factors that result in a
more accurate computations, especially when integrating
with big time steps and for long times. To this end,
we analyze the structural properties of the inelastic flow
governed by a classical material model of finite viscoplasticity.
The material model is
based on the multiplicative decomposition of the
deformation gradient into inelastic and elastic parts.
For simplicity, no hardening behavior is considered in this paper.
However, the proposed methodology can be generalized to
cover more complicated material behavior as well.\footnote{ Using a series of numerical tests, it was shown in
\cite{Shutov1} that the use of geometric integrators allows
to eliminate the error accumulation even in the case of
a more complex material behaviour with nonlinear isotropic and kinematic hardening.
In general, however, the construction of consistent integration procedures
for the finite strain inelasticity is still an open problem (cf. \cite{Vadillo}).}

We pay especial attention to the \emph{exponential stability}
of the inelastic flow, which is the key notion of the current study.
We say that the solution to a Cauchy problem is exponentially stable, if
for small perturbations of initial data,
an exponential decay estimate holds
(see Section 2.1).
From mechanical standpoint, the exponential stability
implies fading memory behavior.\footnote{
As Truesdell and Noll \cite{Truesdell} put it, "Deformations that occurred in the distant past
should have less influence in determining the present stress than those that occurred in the recent past".}
Moreover, the exponential stability is deeply connected to
\emph{contractivity} (B-stability) of the system of equations, which can be used
for stability analysis of numerical algorithms
(see the monograph by Simo and Hughes \cite{SimHug}).

The main conclusions of this paper regarding
the problem of finite viscoplasticity are as follows.
\begin{itemize}
\item The exact solution is \emph{exponentially
stable} with respect to small perturbations of initial data, if
the incompressibility constraint is not violated.
\item In the case of exponential stability, the numerical error
is \emph{uniformly bounded}. In particular, there is
\emph{no error accumulation} even within large time periods.
\item If the incompressibility constraint is violated by some
numerical algorithm, then, in general,
the numerical error tends to \emph{accumulate over time}.
\end{itemize}

There exists a rich mathematical literature dealing with existence,
uniqueness, regularity, and asymptotic behavior of solutions
for certain plascticity/viscoplasticity problems
in the context of infinitesimal strains (see \cite{Alber, HanW, Ionescu} and references therein).
A class of material models of monotone type which
includes the class of generalized standard materials
was defined and analyzed in \cite{Alber}. In the context of finite viscoplasticity, however, only
few theoretical works exist.
Some preliminary investigations have been made by Neff in \cite{Neff}.

In this paper, we analyze the well-known
material model of finite viscoplasticity.
The stability is proved analogously to the classical Lyapunov approach, based on
the use of Lyapunov-candidate-functions.
In fact, the hyperelastic potential is used to construct a suitable Lyapunov candidate
(cf. \cite{SimHug}).

The paper is organized as follows.
In Section 2, we define the notion of exponential stability and prove the main theorem,
which states that the numerical error is uniformly bounded if the exact
solution is exponentially stable. A simple one-dimensional example is presented.
In the next section, a classical material model of finite viscoplasticity
is formulated in the reference configuration.
The change of the reference configuration is likewise discussed.
Section 4 contains the definition and analysis of the distance between two solutions
in terms of energy (Lyapunov candidate). Next, the time-evolution of the distance is evaluated and the
exponential stability of the exact solution is proved.
Finally, the results of numerical tests are presented, which illustrate
the excellent accuracy and
convergence characteristics of geometric integrators.


We conclude this introduction with a few words regarding notation.
Expression $a := b$ means $a$ is defined to be another name for $b$.
Throughout this article, bold-faced symbols denote first- and second-rank tensors in $\mathbb{R}^3$.
A coordinate-free tensor setting is used in this paper (cf. \cite{Itskov, Shutov2}).
The scalar product of two second rank tensors is defined by
$\mathbf{A} : \mathbf{B} = \text{tr}(\mathbf{A} \mathbf{B}^{\text{T}})$.
This scalar product gives rise to the norm by
$\| \mathbf A \| :=  \sqrt{ \mathbf{A} : \mathbf{A}}$.
Moreover, we denote by $\| \cdot \|^{*} $ the induced norm of a tensor
\begin{equation}\label{Openo}
\| \mathbf A \|^{*} := \max_{\|\mathbf x\|_2=1} \| \mathbf A \mathbf x \|_2, \quad \|\mathbf x\|_2 := \sqrt{\mathbf x \cdot \mathbf x}.
\end{equation}
The overline $\overline{(\cdot)}$ stands for the unimodular part of a tensor
\begin{equation}\label{UniPro}
\overline{\mathbf{A}}:=(\det \mathbf{A})^{-1/3} \mathbf{A}.
\end{equation}
The deviatoric part of a tensor is defined as
$\mathbf A^{\text{D}} := \mathbf A - \frac{1}{3} \text{tr}(\mathbf A) \mathbf 1$.
The notation $O$ stands for "Big-O" Landau symbol:
$f(x)= O (g(x)) \ \text{as} \ x \ \rightarrow x_0$ iff there exists  $C < \infty \ \text{such that} \
\|f(x)\| \leq C \|g(x)\| \ \text{as} \ x \ \rightarrow x_0$.
The inequality $f(x) \leq O (g(x))$ is understood as follows:
there exists $\acute{f}(x) = O (g(x))$ such that $f(x) \leq \acute{f}(x)$.

\section{Differential equations on manifolds and exponential stability}

\subsection{General definitions}

Let us consider the Cauchy problem for a smooth function $y(t) \in \mathbb{R}^n$
\begin{equation}{\label{Cauchy}}
\dot{y} (t) = f(y(t), d(t)), \quad y(t_0)= y_0.
\end{equation}
Here, the initial value $y_0$ and the function $d(t)$ are
supposed to be given.\footnote{The system \eqref{Cauchy} is
a system with input, and $d(t)$ is interpreted as a forcing function.}
Denote the exact solution to \eqref{Cauchy} by $\tilde{y}(t,y_0,t_0)$. In particular, we have
\begin{equation}{\label{note1}}
\tilde{y}(t_0,y_0,t_0) = y_0.
\end{equation}
Suppose that all solutions lie on some manifold $M \subset \mathbb{R}^n$
\begin{equation}{\label{manifold}}
\tilde{y}(t,y_0,t_0) \in M, \quad \text{for all} \ t \geq t_0, \ y_0 \in M.
\end{equation}
Then we say that \eqref{Cauchy} is a \emph{differential equation on
the manifold} $M$ (cf. \cite{Hair, Hair2}).

Next, we say that the solution $y(t)$ to the problem \eqref{Cauchy} is \emph{locally
exponentially stable} on $M$, if there exist
$\delta >0, \ \gamma > 0, \ C_1 < \infty$, such that the following decay estimate holds
\begin{equation}{\label{decay}}
\| \tilde{y}(t,y^{(1)}_0,t_0) - \tilde{y}(t,y^{(2)}_0,t_0) \|  \leq  C_1 \ \e^{\displaystyle - \gamma (t-t_0)} \ \| y^{(1)}_0 - y^{(2)}_0 \|,
\end{equation}
for all $t_0 \geq 0$, $y^{(1)}_0, y^{(2)}_0 \in M$ such that $\| y^{(1)}_0 - y(t_0) \| \leq \delta$,
$\| y^{(2)}_0 - y(t_0) \| \leq \delta$.

We note that somewhat different interpretation of the exponential stability
can be met in the literature as well (cf., for example, Section 2.5 of \cite{Kato}).

Next, let us consider a numerical algorithm which solves \eqref{Cauchy} on the
time interval $[0, T]$. Denote by ${}^n y$ the numerical solutions
at time instances ${}^n t$, where
$0={}^0 t < {}^1 t < {}^2 t < ... < {}^N t = T$, and ${}^0 y = y_0$. Suppose that the error on the step is
bounded by the second power of the step size. More precisely
\begin{equation}{\label{ErrStep}}
\| \tilde{y}({}^{n+1} t,{}^n y,{}^n t) - {}^{n+1} y \| \leq C_2 ({}^{n+1} t - {}^n t)^2,
\end{equation}
where $C_2 < \infty$ (cf. figure \ref{fig1}).
For simplicity, we will consider
constant time-steps only: $\Delta t = {}^{n+1} t - {}^n t = const$.

\subsection{Main theorem}

With definitions from previous section we formulate the following theorem.

\textbf{Theorem 1.}

Let $y(t) = \tilde{y}(t,y_0, 0)$ be the exact solution. Suppose that
conditions \eqref{decay} and \eqref{ErrStep} hold.
Moreover, suppose that the numerical solution of problem \eqref{Cauchy} lies
exactly on $M$. Then there exist a constant $C < \infty$ such that
\begin{equation}{\label{theorem}}
\| {}^{n} y - y({}^n t) \| \leq C \Delta t, \quad \text{as} \ \Delta t \rightarrow 0.
\end{equation}
Here, the constant $C$ does not depend on the size of the time interval $[0, T]$.

\emph{Proof.}
The proof is a modification of the standard error analysis (cf. \cite{Bahvalov}).
In this paper we prove the theorem under assumption that $\delta = \infty$.
The proof can be easily generalized to cover arbitrary values of $\delta > 0$
by using mathematical induction and by assuming $\Delta t \leq \gamma \delta / (2 \max ({C_1,1}) \ C_2)$.

\begin{figure}\centering
\psfrag{A0}[m][][1][0]{${}^0 t$}
\psfrag{A1}[m][][1][0]{${}^1 t$}
\psfrag{A2}[m][][1][0]{${}^2 t$}
\psfrag{A3}[m][][1][0]{${}^3 t$}
\psfrag{A4}[m][][1][0]{${}^4 t$}
\psfrag{B}[m][][1][0]{${}^4 y$}
\psfrag{C}[m][][1][0]{$\tilde{y}({}^4 t,{}^3 y,{}^3 t)$}
\psfrag{D}[m][][1][0]{$\tilde{y}({}^4 t,{}^2 y,{}^2 t)$}
\psfrag{E}[m][][1][0]{$\tilde{y}({}^4 t,{}^1 y,{}^1 t)$}
\psfrag{F}[m][][1][0]{$y ({}^4 t)$}
\psfrag{G}[m][][1][0]{Exact solution}
\psfrag{H}[m][][1][0]{Numerical solution}
\scalebox{0.9}{\includegraphics{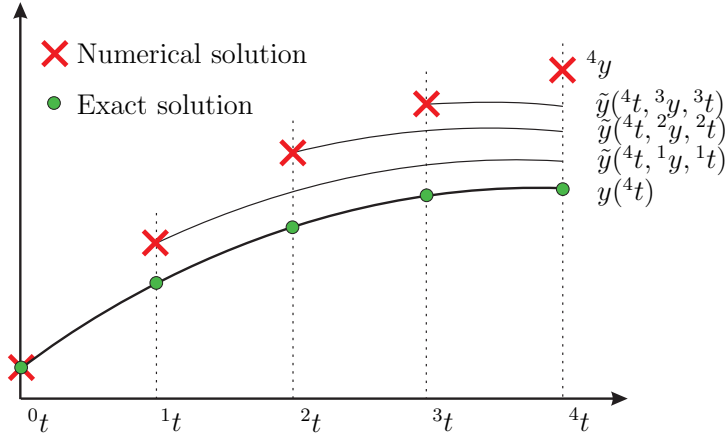}}
\caption{Analysis of error accumulation. \label{fig1}}
\end{figure}

First, note that 
$\tilde{y}({}^n t,{}^{0} y,{}^{0} t) = y({}^n t)$. Thus, (cf. figure \ref{fig1})
\begin{multline}{\label{proof1}}
\| {}^{n} y - y({}^n t) \| \leq
\| {}^{n} y - \tilde{y}({}^n t,{}^{n-1} y,{}^{n-1} t) \| +
\| \tilde{y}({}^n t,{}^{n-1} y,{}^{n-1} t) - \tilde{y}({}^n t,{}^{n-2} y,{}^{n-2} t) \| + ... \\
+ \| \tilde{y}({}^n t,{}^{1} y,{}^{1} t) -  \tilde{y}({}^n t,{}^{0} y,{}^{0} t) \|.
\end{multline}
Next, from \eqref{decay} we obtain for all $k = 1,2,..., n-1$
\begin{equation}{\label{decay2}}
\| \tilde{y}({}^n t,{}^{k} y,{}^{k} t) - \tilde{y}({}^n t,{}^{k-1} y,{}^{k-1} t) \|
\leq  C_1 \ \e^{\displaystyle - \gamma ({}^n t - {}^{k} t)} \
\| \tilde{y}({}^k t,{}^{k} y,{}^{k} t) - \tilde{y}({}^k t,{}^{k-1} y,{}^{k-1} t) \|.
\end{equation}
Substituting \eqref{decay2} in \eqref{proof1}, we get
\begin{multline}{\label{proof3}}
\| {}^{n} y - y({}^n t) \| \leq \\
\| {}^{n} y - \tilde{y}({}^n t,{}^{n-1} y,{}^{n-1} t) \| + C_1
\sum_{k=1}^{n-1} \e^{\displaystyle - \gamma ({}^n t - {}^{k} t)} \
\| \tilde{y}({}^k t,{}^{k} y,{}^{k} t) - \tilde{y}({}^k t,{}^{k-1} y,{}^{k-1} t) \|.
\end{multline}
Obviously, $\tilde{y}({}^k t,{}^{k} y,{}^{k} t) = {}^{k} y$.
Without loss of generality, we can assume that $C_1 \geq 1$.
Next, substituting error
estimation \eqref{ErrStep} into \eqref{proof3}, we get
\begin{equation}{\label{proof5}}
\| {}^{n} y - y({}^n t) \| \leq   C_1 \ C_2 \ (\Delta t)^2
\sum_{k=1}^{n} \e^{\displaystyle - \gamma ({}^n t - {}^{k} t)}.
\end{equation}
But, ${}^n t - {}^{k} t = (n-k) \Delta t$. Thus, taking into account the well-known expression
for an infinite geometric series ($\sum_{i=0}^{\infty} r^i = 1/(1-r)$ for $|r | < 1 $ ),
we get for small $\Delta t$
\begin{equation}{\label{proof6}}
\sum_{k=1}^{n} \e^{\displaystyle - \gamma ({}^n t - {}^{k} t)} \leq
\sum_{i=0}^{\infty} \e^{\displaystyle - i \gamma \Delta t} =
\frac{\displaystyle 1}{\displaystyle 1 - \e^{ - \gamma \Delta t}}
= \frac{\displaystyle 1}{\displaystyle \gamma \Delta t + O((\Delta t)^2 )} \leq 2 \frac{\displaystyle 1}{\displaystyle \gamma \Delta t}.
\end{equation}
Finally, it follows from \eqref{proof5}, \eqref{proof6}
\begin{equation}{\label{proof7}}
\| {}^{n} y - y({}^n t) \| \leq  \frac{\displaystyle 2 C_1 \ C_2}{\displaystyle \gamma}
\ \Delta t, \quad \text{as} \ \Delta t \rightarrow 0 \quad \quad \quad \quad \quad \quad  \blacksquare
\end{equation}

\textbf{Remark 1.}
The proof is essentially based on the assumption that ${}^{k} y \in M$.
In general, if the numerical solution ${}^{k} y$ leaves the manifold $M$,
the decay estimation \eqref{decay} is not valid.

\textbf{Remark 2.}
The theorem states that the error is uniformly bounded in the
case of exponential stability. Thus, there is no error accumulation in the sense that
the constant $C$ in \eqref{theorem} does not depend on the overall time $T$.
Moreover, let $\epsilon > 0$ be some small value.
By choosing  $\Delta t \leq \gamma \epsilon / (2 C_1 C_2)$
the numerical error $\| {}^{n} y - y({}^n t) \|$ is guaranteed to be less than $\epsilon$.

\textbf{Remark 3.}
If the exponential stability is replaced by the
assumption that the right-hand side of \eqref{Cauchy}
is a smooth function of $y$, a weaker error estimation is valid (cf. \cite{Bahvalov})
\begin{equation}{\label{weaker}}
\| {}^{n} y - y({}^n t) \| \leq  C \e^{L T} T
\ \Delta t, \quad \text{as} \ \Delta t \rightarrow 0,
\end{equation}
where $L = \sup  \| f_y \|$.
The effect of growing multiplier on the right hand side of \eqref{weaker}
is referenced to as an effect of error accumulation.
In that case, in order to guaranty
a sufficient accuracy,
the upper bound for $\Delta t$ must depend on $T$. That makes the practical
solution of some problems extremely expensive for large values of $T$.

\subsection{One-dimensional example}

\begin{figure}
\psfrag{A}[m][][1][0]{$\varepsilon_{\text{i}}$}
\psfrag{B}[m][][1][0]{$\varepsilon_{\text{e}}$}
\psfrag{C}[m][][1][0]{$\varepsilon$}
\psfrag{B1}[m][][1][0]{$\varepsilon_{\text{i}}$}
\psfrag{A1}[m][][1][0]{$t$}
\psfrag{E}[m][][1][0]{(a)}
\psfrag{F}[m][][1][0]{(b)}
\scalebox{1.0}{\includegraphics{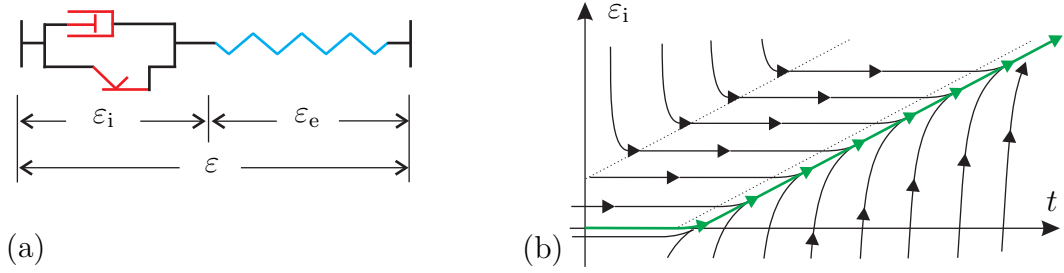}}
\caption{Rheological model (a), and inelastic flow
under monotonic loading (b). \label{fig2}}
\end{figure}

Let us consider a simple example which illustrates the notion of exponential stability.
We examine the response of a one-dimensional viscoplastic device shown in Figure \ref{fig2} (a).
The closed system of (constitutive) equations is as follows:

The total strain is decomposed into elastic part $\varepsilon_{\text{e}}$, and
inelastic part $\varepsilon_{\text{i}}$
\begin{equation}{\label{examp1}}
\varepsilon = \varepsilon_{\text{e}} + \varepsilon_{\text{i}}.
\end{equation}
The stress $\sigma$ on the elastic spring is governed by elasticity law ($E > 0$).
\begin{equation}{\label{examp2}}
\sigma = E \varepsilon_{\text{e}}.
\end{equation}
The time derivative of the inelastic strain is given by
\begin{equation}{\label{examp3}}
\dot{\varepsilon}_{\text{i}} = \frac{1}{\eta} \langle  f  \rangle \frac{\sigma}{|\sigma|}, \quad
f: = |\sigma| - K, \quad \langle  x  \rangle := \max (x,0),
\end{equation}
where material constants $K > 0$ and
$\eta > 0$ are referred to as yield stress and viscosity, respectively.

In order to use the results of previous subsections, we rewrite the problem
in the form 
\begin{equation}{\label{examp5}}
\dot{\varepsilon}_{\text{i}} = \dot{\varepsilon}_{\text{i}} (\varepsilon_{\text{i}}, \varepsilon(t)) = \frac{1}{\eta}
\langle E |\varepsilon(t) - \varepsilon_{\text{i}}| - K \rangle \ \text{sign}(\varepsilon(t) - \varepsilon_{\text{i}} ).
\end{equation}
Let $\varepsilon^{(1)}_{\text{i}} (t)$ and
$\varepsilon^{(2)}_{\text{i}} (t)$ be to two solutions to \eqref{examp5}.
Following \cite{SimHug}, we recall that \\
$\sqrt{\frac{1}{2} E (\varepsilon^{(1)}_{\text{i}} -
\varepsilon^{(2)}_{\text{i}})^2}$
defines an \emph{energy norm} which is the natural norm
for the problem under consideration.\footnote
{It is known (see \cite{SimHug}) that
$\frac{1}{2} E (\varepsilon^{(1)}_{\text{i}}(t)- \varepsilon^{(2)}_{\text{i}}(t))^2$
is not increasing. This effect is referenced to as \emph{contractivity}.}
Next, we consider a monotonic loading
\begin{equation}{\label{examp4}}
\varepsilon (t) = \dot{\varepsilon} t, \quad \dot{\varepsilon} = const >0.
\end{equation}
Let us show that the exact solution satisfying the initial
condition $\varepsilon_{\text{i}}=0 $ is
exponentially stable.\footnote{For the current example, the geometric property $y \in M$
is trivial: we put $M = \mathbb{R}.$}
Without loss of generality, we can assume that $t_0 = 0$ in estimation \eqref{decay}.
If $|\varepsilon^{(k)}_{\text{i}}(0) - 0| \leq \delta$ for $\text{k} \in \{1,2\}$, then there exists time instance $t' = t' (\delta)$
such that the condition $f \geq f_0 > 0$ holds for both solutions ($\varepsilon^{(1)}_{\text{i}} (t)$ and
$\varepsilon^{(2)}_{\text{i}} (t)$), if $t \geq t'$ (see Figure \ref{fig2} (b)).
Then, under that assumption
\begin{equation}{\label{examp6}}
\frac{\partial \dot{\varepsilon}_{\text{i}}
(\varepsilon_{\text{i}}, \varepsilon(t))}{\partial  \varepsilon_{\text{i}}} = - \frac{E}{\eta}, \quad
\dot{\varepsilon}_{\text{i}}
(\varepsilon^{(1)}_{\text{i}}, \varepsilon(t)) - \dot{\varepsilon}_{\text{i}}
(\varepsilon^{(2)}_{\text{i}}, \varepsilon(t)) = - \frac{E}{\eta} (\varepsilon^{(1)}_{\text{i}} - \varepsilon^{(2)}_{\text{i}}).
\end{equation}
Therefore, we get from \eqref{examp6}
\begin{equation}{\label{examp7}}
\Big( \frac{1}{2} E (\varepsilon^{(1)}_{\text{i}} - \varepsilon^{(2)}_{\text{i}})^2 \Big)^{\cdot} =
E (\varepsilon^{(1)}_{\text{i}} - \varepsilon^{(2)}_{\text{i}})
(\dot{\varepsilon}^{(1)}_{\text{i}} - \dot{\varepsilon}^{(2)}_{\text{i}}) =
-\frac{E^2}{\eta} (\varepsilon^{(1)}_{\text{i}} - \varepsilon^{(2)}_{\text{i}})^2, \ \text{for} \ t \geq t'.
\end{equation}

Due to the contractivity (for details see \cite{SimHug}),
$\frac{1}{2} E (\varepsilon^{(1)}_{\text{i}}(t') - \varepsilon^{(2)}_{\text{i}}(t'))^2  \leq
 \frac{1}{2} E (\varepsilon^{(1)}_{\text{i}}(0) - \varepsilon^{(2)}_{\text{i}}(0))^2 $. Moreover, integrating
\eqref{examp7} over $[t', t]$ and taking the contractivity into account, we get
\begin{equation}{\label{examp8}}
\frac{1}{2} E (\varepsilon^{(1)}_{\text{i}}(t) - \varepsilon^{(2)}_{\text{i}}(t))^2 \leq
\frac{1}{2} E (\varepsilon^{(1)}_{\text{i}}(0) - \varepsilon^{(2)}_{\text{i}}(0))^2 \e^{-\frac{2 E}{\eta} (t- t')}.
\end{equation}
Taking the square root of both sides we obtain the required exponential decay estimation \eqref{decay} with
$C_1 = \e^{\frac{\displaystyle E t'}{\displaystyle \eta}} < \infty$ and $\gamma = \frac{\displaystyle E}{\displaystyle \eta} > 0$.


\section{Material model of multiplicative viscoplasticity}
Let us consider a classical material model of finite viscoplasticity (see, for example, \cite{Haupt}).

\subsection{Constitutive equations}
The model is based on the multiplicative split of the deformation gradient $\mathbf F$
\begin{equation}\label{split1}
\mathbf F = \hat{\mathbf F}_{\text{e}} \mathbf F_{\text{i}}.
\end{equation}
Here, $\hat{\mathbf F}_{\text{e}}$ and $\mathbf F_{\text{i}}$
stand for \emph{elastic} and \emph{inelastic} parts, respectively (see \cite{Kroener, Lee}).
The multiplicative split can be motivated by
the idea of a local elastic unloading. A somewhat more
consistent motivation
can be derived from the
concept of material isomorphism \cite{Bertram}.

Along with the well-known
right Cauchy-Green tensor
${\mathbf C}=\mathbf F^{\text{T}} \mathbf F$, we introduce
a strain-like internal variable (inelastic right Cauchy-Green tensor) as
\begin{equation}\label{intIntv}
{\mathbf C}_{\text{i}}=\mathbf F^{\text{T}}_{\text{i}} \mathbf F_{\text{i}}.
\end{equation}

In this paper we consider strain-driven processes. More precisely, we assume
the deformation history ${\mathbf C}(t)$ to be given.
The material response in the time interval $t \in [0,T]$ is governed by
the following ordinary differential equation with initial condition
\begin{equation}\label{prob1}
\dot{\mathbf C}_{\text{i}} = 2 \frac{\displaystyle
\lambda_{\text{i}}}{\displaystyle \mathfrak{F}}
 \big( \mathbf C \tilde{\mathbf T} \big)^{\text{D}} \mathbf C_{\text{i}}, \quad
\mathbf C_{\text{i}}|_{t=0} = \mathbf C_{\text{i}}^0, \
\det \mathbf C_{\text{i}}^0 =1, \ \mathbf C_{\text{i}}^0 \in Sym.
\end{equation}
Here, the 2nd Piola-Kirchhoff tensor $\tilde{\mathbf T}$,
the norm of the driving force $\mathfrak{F}$, and the inelastic multiplier $\lambda_{\text{i}}$
are functions of $(\mathbf C, \mathbf C_{\text{i}})$, given by
\begin{equation}\label{prob2}
\tilde{\mathbf T} =
2 \rho_{\scriptscriptstyle \text{R}}
\frac{\displaystyle \partial \psi_{\text{el}}
(\mathbf C {\mathbf C_{\text{i}}}^{-1})}
{\displaystyle \partial
\mathbf{C}}\big|_{\mathbf C_{\text{i}} = \text{const}},
\end{equation}
\begin{equation}\label{prob3}
\lambda_{\text{i}}= \frac{\displaystyle 1}{\displaystyle
\eta}\Big\langle \frac{\displaystyle 1}{\displaystyle k_0}
f \Big\rangle^{m}, \quad
f:= \mathfrak{F}- \sqrt{\frac{2}{3}} K, \quad
\mathfrak{F}:= \sqrt{\text{tr}
\big[ \big( \mathbf C \tilde{\mathbf T} \big)^{\text{D}} \big]^2 }.
\end{equation}
The material parameters
$\rho_{\scriptscriptstyle \text{R}} > 0$,
$\eta \geq 0$, $m \geq 1$, $K > 0$, and
the isotropic real-valued function $\psi_{\text{el}}$ are assumed to be known;
$k_0 >0$ is used to get a dimensionless term in the bracket.

\textbf{Remark.}
The right Cauchy strain tensor $\mathbf C$ is symmetric.
Since the function $\psi_{\text{el}}$ is isotropic,
it makes no difference whether the derivative in \eqref{prob2} is understood
as a general derivative or as a derivatives with respect
to a symmetric tensor (cf. \cite{Shutov2}).

Next, we remark that the right-hand side in $\eqref{prob1}_1$
is symmetric (cf. \cite{Shutov2}).
Moreover, taking into account the property
$\text{tr}(\mathbf{A} \mathbf{B}) = \text{tr}(\mathbf{B} \mathbf{A})$ and combining the Jacobi
formulae $\frac{\displaystyle \partial \text{det}(\mathbf{A})}{\displaystyle
\partial  \mathbf{A}} = \text{det}(\mathbf{A}) \ \mathbf{A}^{-\text{T}}$
with the evolution equation $\eqref{prob1}_1$, we get
\begin{equation}\label{incompr1}
\big(\det \mathbf C_{\text{i}} \big)^{\cdot} = \frac{\displaystyle \partial
\text{det}(\mathbf{\mathbf C_{\text{i}}})}{\displaystyle
\partial  \mathbf{\mathbf C_{\text{i}}}} : \dot{\mathbf C}_{\text{i}}  =
2 \frac{\displaystyle \lambda_{\text{i}}}{\displaystyle \mathfrak{F}}
\text{det}(\mathbf{\mathbf C_{\text{i}}}) \ \mathbf C_{\text{i}}^{-1} :
 \big( \mathbf C \tilde{\mathbf T} \big)^{\text{D}} \mathbf C_{\text{i}} =
 2 \frac{\displaystyle
\lambda_{\text{i}}}{\displaystyle \mathfrak{F}}
\text{tr}\big[\big( \mathbf C \tilde{\mathbf T} \big)^{\text{D}}\big] = 0.
\end{equation}

Therefore, the exact solution of
\eqref{prob1} -- \eqref{prob3} has the following geometric property
\begin{equation}\label{geopr0}
\mathbf{C}_{\text{i}} \in M, \quad
M := \big\{ \mathbf B \in Sym \ | \ \text{det} \mathbf B =1 \big\}.
\end{equation}
We note that the current material model is \emph{thermodynamically consistent}.
That means that the Clausius-Duhem inequality holds for arbitrary
mechanical loadings
\begin{equation}\label{cld}
\delta_{\text{i}} := \frac{1}{2 \rho_{\scriptscriptstyle \text{R}}}
\tilde{\mathbf T} : \dot{\mathbf C} -
\big( \psi_{\text{el}} (\mathbf C {\mathbf C_{\text{i}}}^{-1}) \big)^{\cdot} \geq 0.
\end{equation}
In particular, we get a reduced inequality for relaxation processes ($\mathbf C = \text{const}$)
\begin{equation}\label{decayCld}
\big( \psi_{\text{el}} (\mathbf C {\mathbf C_{\text{i}}}^{-1}) \big)^{\cdot} \ \leq 0.
\end{equation}
One mathematical interpretation of this inequality will be discussed in Section 4.1.

To be definite, we use the following expression for the free energy density $\psi_{\text{el}}$
(generalized Neo-Hooke model \cite{Haupt})
\begin{equation}\label{spec1}
\rho_{\scriptscriptstyle \text{R}}  \psi_{\text{el}}(\mathbf{A}):=
\frac{k}{2}\big(\text{ln}\sqrt{\text{det} \mathbf{A}} \big)^2+
\frac{\mu}{2} \big( \text{tr} \overline{\mathbf{A}} - 3 \big),
\end{equation}
where $k >0$, $\mu > 0$ are known material constants (bulk modulus and shear modulus, respectively).

Substituting \eqref{spec1} in \eqref{prob2} we get the 2nd Piola-Kirchhoff stress tensor in the form
\begin{equation}\label{spec2}
\tilde{\mathbf T} = k \ \text{ln}\sqrt{\text{det} (\mathbf C)} \
  \mathbf C^{-1} + \mu \ \mathbf C^{-1} (\overline{\mathbf C}
  \mathbf C_{\text{i}}^{-1})^{\text{D}}.
\end{equation}

In what follows we analyze the exponential stability of the exact solution $\mathbf{C}_{\text{i}}(t)$.

\subsection{Change of reference configuration}

In order to simplify the analysis of the material model, we may need to
rewrite the constitutive equation with respect to some "new" local
reference configuration $\mathbf F_{0}$. In what follows, we suppose that this configuration
is isochoric, i.e. $\det (\mathbf F_{0})=1$.
The "new" deformation gradient, right Cauchy tensor, and inelastic right Caushy
tensor are given by
\begin{equation}\label{ConfigChange1}
\mathbf F^{\text{new}} := \mathbf F \mathbf F^{-1}_{0}, \quad
\mathbf C^{\text{new}} := \mathbf F^{-\text{T}}_{0} \mathbf C \mathbf F^{-1}_{0}, \quad
\mathbf C^{\text{new}}_{\text{i}} := \mathbf F^{-\text{T}}_{0} \mathbf C_{\text{i}} \mathbf F^{-1}_{0}.
\end{equation}
The 2nd Piola-Kirchhoff tensor $\tilde{\mathbf T}$,
the norm of the driving force $\mathfrak{F}$, the inelastic multiplier $\lambda_{\text{i}}$,
and the overstress $f$ are transformed as follows
\begin{equation}\label{ConfigChange3}
\tilde{\mathbf T}^{\text{new}} := \mathbf F_{0} \tilde{\mathbf T} \mathbf F^{\text{T}}_{0}, \quad
\mathfrak{F}^{\text{new}} := \mathfrak{F}, \quad \lambda_{\text{i}}^{\text{new}} := \lambda_{\text{i}}
\quad f^{\text{new}} := f.
\end{equation}
Since $\psi_{\text{el}}$ is isotropic, $\psi_{\text{el}} (\mathbf{A}
\mathbf{B}) = \psi_{\text{el}} (\mathbf{B} \mathbf{A})$.
Using that property, it can be checked that $\psi_{\text{el}} (\mathbf C {\mathbf C_{\text{i}}}^{-1})$ is
invariant under the change of reference configuration
\begin{equation}\label{ConfigChange2}
\psi_{\text{el}} (\mathbf C^{\text{new}}
(\mathbf C^{\text{new}}_{\text{i}})^{-1}) =
\psi_{\text{el}} (\mathbf C {\mathbf C_{\text{i}}}^{-1}).
\end{equation}
The closed system of equations with respect to the new reference configuration is
obtained from \eqref{prob1} --- \eqref{prob3} by replacing all quantities by
their "new" counterparts.

\section{Analysis of exponential stability for multiplicative viscoplasticity}

\subsection{Measuring the distance between solutions in terms of energy}

Suppose that $\mathbf{C}^{(1)}_{\text{i}}(t)$ and $\mathbf{C}^{(2)}_{\text{i}}(t)$
are two solutions to the problem \eqref{prob1} --- \eqref{prob3} (with the
same forcing function $\mathbf{C}(t)$).
Next, suppose that there exists a constant $L < \infty$ such that
\begin{equation}\label{BoundL}
\| (\mathbf{C}^{(k)}_{\text{i}})^{1/2} (t) \| < L, \quad \| (\mathbf{C}^{(k)}_{\text{i}})^{-1/2} (t) \| < L \quad \text{for all} \ t> 0, \
k \in \{1,2\}.
\end{equation}
We introduce the following measure of distance between two solutions in terms of energy
\footnote{The relation \eqref{EnerDist1} can be seen as a generalization of the energy norm
$\sqrt{\frac{1}{2} E (\varepsilon^{(1)}_{\text{i}} -
\varepsilon^{(2)}_{\text{i}})^2}$ (cf. Section 2.3).}

\begin{equation}\label{EnerDist1}
\text{Dist}(\mathbf{C}^{(1)}_{\text{i}}, \mathbf{C}^{(2)}_{\text{i}})
 := \sqrt{ \rho_{\scriptscriptstyle \text{R}} \psi_{\text{el}} \big(\mathbf{C}^{(1)}_{\text{i}}
(\mathbf{C}^{(2)}_{\text{i}})^{-1}\big)}.
\end{equation}
This measure has the following properties:

(i) Invariance under the change of reference configuration
\begin{equation}\label{EnerDist2}
\text{Dist}\big( (\mathbf{C}^{(1)}_{\text{i}})^{\text{new}},
(\mathbf{C}^{(2)}_{\text{i}})^{\text{new}} \big)
 = \text{Dist}\big(\mathbf{C}^{(1)}_{\text{i}}, \mathbf{C}^{(2)}_{\text{i}}\big).
\end{equation}

(ii) For small $\mathbf{C}^{(1)}_{\text{i}} - \mathbf{C}^{(2)}_{\text{i}}$,
there exist constants $C_3 > 0$ and $C_4 < \infty$ such that
\begin{equation}\label{EnerDist3}
C_3 \ \|\mathbf{C}^{(1)}_{\text{i}} - \mathbf{C}^{(2)}_{\text{i}} \|
\leq \text{Dist}(\mathbf{C}^{(1)}_{\text{i}}, \mathbf{C}^{(2)}_{\text{i}}) \leq C_4
\|\mathbf{C}^{(1)}_{\text{i}} - \mathbf{C}^{(2)}_{\text{i}} \|.
\end{equation}

(iii) For all $\mathbf{C}^{(1)}_{\text{i}}(t), \mathbf{C}^{(2)}_{\text{i}}(t) \in M$
we have $\text{Dist}(\mathbf{C}^{(1)}_{\text{i}},
\mathbf{C}^{(2)}_{\text{i}}) \geq 0$ and
\begin{equation}\label{EnerDist4}
\text{Dist}(\mathbf{C}^{(1)}_{\text{i}},
\mathbf{C}^{(2)}_{\text{i}}) = 0, \ \text{if and only if} \ \ \mathbf{C}^{(1)}_{\text{i}} = \mathbf{C}^{(2)}_{\text{i}}.
\end{equation}

\textbf{\emph{Proof.}}

(i): Identity \eqref{EnerDist2} can be proved similarly to the invariance property \eqref{ConfigChange2}.

(ii): First, it follows from \eqref{spec1} that for small $\mathbf{\Delta}$ we have (see Appendix A)
\begin{equation}\label{EnerDist5}
\rho_{\scriptscriptstyle \text{R}} \psi_{\text{el}} (\mathbf{1} + \mathbf{\Delta}) =
\frac{k}{8} (\text{tr} \mathbf{\Delta})^2 +
\frac{\mu}{4} \text{tr} \big( (\mathbf{\Delta}^{\text{D}} )^2 \big) + O( \Delta^3 ),
\end{equation}
where $\Delta = \| \mathbf{\Delta} \|$. Note that
$\text{tr} \big( (\mathbf{\Delta}^{\text{D}} )^2 \big) = \|\mathbf{\Delta}^{\text{D}}\|^2$
for $\mathbf{\Delta} \in Sym$. Thus, $\frac{k}{8} (\text{tr} \mathbf{\Delta})^2 +
\frac{\mu}{4} \text{tr} \big( (\mathbf{\Delta}^{\text{D}} )^2 \big)$ is a norm on $Sym$.
Since all norms on $Sym$ are equivalent,
there exist constants $C'_3 > 0$, $C'_4 < \infty$ such that
for small $\mathbf{\Delta} \in Sym$ we have
\begin{equation}\label{EnerDist6}
C'_3 \| \mathbf{\Delta} \|^2 \leq \rho_{\scriptscriptstyle \text{R}} \psi_{\text{el}}
(\mathbf{1} + \mathbf{\Delta}) \leq C'_4 \| \mathbf{\Delta} \|^2.
\end{equation}
Next, due to the property
\begin{equation}\label{EnerDist62}
\psi_{\text{el}} (\mathbf{A}
\mathbf{B}) = \psi_{\text{el}} (\mathbf{B} \mathbf{A}),
\end{equation}
we have
\begin{equation}\label{EnerDist7}
\text{Dist}(\mathbf{C}^{(1)}_{\text{i}}, \mathbf{C}^{(2)}_{\text{i}}) =
\sqrt{\rho_{\scriptscriptstyle \text{R}} \psi_{\text{el}} \big( (\mathbf{C}^{(2)}_{\text{i}})^{-1/2} \
\mathbf{C}^{(1)}_{\text{i}} \
(\mathbf{C}^{(2)}_{\text{i}})^{-1/2} \big)}.
\end{equation}
Moreover, taking into account that
$\| \mathbf{A} \mathbf{B} \|_* \leq \| \mathbf{A}  \|_*  \  \| \mathbf{B} \|_*$,
and that the norms $\| \cdot \|_*$ and $\| \cdot \|$ are equivalent, we get
\begin{equation}\label{EnerDist8}
\| \mathbf{A} \mathbf{B} \mathbf{C}\| \leq \acute{C} \| \mathbf{A}\| \ \| \mathbf{B}\| \ \| \mathbf{C}\|,
\end{equation}
with some constant $\acute{C} < \infty$. Thus,
\begin{multline}\label{EnerDist9}
\big \| (\mathbf{C}^{(2)}_{\text{i}})^{-1/2} \
\mathbf{C}^{(1)}_{\text{i}} \
(\mathbf{C}^{(2)}_{\text{i}})^{-1/2} - \mathbf{1} \big  \| =
\big \| (\mathbf{C}^{(2)}_{\text{i}})^{-1/2} \
\big( \mathbf{C}^{(1)}_{\text{i}} -  \mathbf{C}^{(2)}_{\text{i}}  \big) \
(\mathbf{C}^{(2)}_{\text{i}})^{-1/2} \big \|  \\ \stackrel{\eqref{EnerDist8}}{\leq}
\acute{C}
\big \| (\mathbf{C}^{(2)}_{\text{i}})^{-1/2} \big \|^2 \
\big \| \mathbf{C}^{(1)}_{\text{i}} -  \mathbf{C}^{(2)}_{\text{i}} \big \|,
\end{multline}
\begin{multline}\label{EnerDist10}
\big \| \mathbf{C}^{(1)}_{\text{i}} -  \mathbf{C}^{(2)}_{\text{i}} \big \| =
\big \| (\mathbf{C}^{(2)}_{\text{i}})^{1/2} \big((\mathbf{C}^{(2)}_{\text{i}})^{-1/2} \
\mathbf{C}^{(1)}_{\text{i}} \
(\mathbf{C}^{(2)}_{\text{i}})^{-1/2} - \mathbf{1}\big)  (\mathbf{C}^{(2)}_{\text{i}})^{1/2}   \big \| \\
\stackrel{\eqref{EnerDist8}}{\leq}
\acute{C} \big \| (\mathbf{C}^{(2)}_{\text{i}})^{1/2} \big \|^2 \
\big \| (\mathbf{C}^{(2)}_{\text{i}})^{-1/2} \
\mathbf{C}^{(1)}_{\text{i}} \
(\mathbf{C}^{(2)}_{\text{i}})^{-1/2} - \mathbf{1} \big \| .
\end{multline}
Further,
substituting $\mathbf{\Delta} = (\mathbf{C}^{(2)}_{\text{i}})^{-1/2} \
\mathbf{C}^{(1)}_{\text{i}} \
(\mathbf{C}^{(2)}_{\text{i}})^{-1/2} - \mathbf{1}$ in
\eqref{EnerDist6}, and combining it with
\eqref{EnerDist7}, \eqref{EnerDist9}, and \eqref{EnerDist10} we get
\begin{equation}\label{EnerDist32}
\tilde{C}_3 \| (\mathbf{C}^{(2)}_{\text{i}})^{1/2} \|^{-2}
\ \|\mathbf{C}^{(1)}_{\text{i}} - \mathbf{C}^{(2)}_{\text{i}} \|
\leq \text{Dist}(\mathbf{C}^{(1)}_{\text{i}}, \mathbf{C}^{(2)}_{\text{i}}) \leq \tilde{C}_4
\| (\mathbf{C}^{(2)}_{\text{i}})^{-1/2} \|^2 \
\|\mathbf{C}^{(1)}_{\text{i}} - \mathbf{C}^{(2)}_{\text{i}} \|.
\end{equation}
Finally, combining \eqref{EnerDist32} with \eqref{BoundL} we obtain \eqref{EnerDist3}.

(iii): We note that $\psi_{\text{el}} (\mathbf{A}) \geq 0$. Moreover,
$\psi_{\text{el}}(\mathbf{A}) = 0$ if and only if $\mathbf{A}=\mathbf{1}$ $\blacksquare$

In view of properties (i) --- (iii), the function Dist is
a natural measure of distance for the problem under consideration.\footnote
{The function Dist is not symmetric: $\text{Dist}(\mathbf{A}, \mathbf{B}) \neq \text{Dist}(\mathbf{B}, \mathbf{A})$.
Symmetrized functions can be defined by
${\text{Dist}}^{Sym}_1(\mathbf{A}, \mathbf{B}) := 1/2 (\text{Dist}(\mathbf{A}, \mathbf{B}) + \text{Dist}(\mathbf{B}, \mathbf{A}))$,
${\text{Dist}}^{Sym}_2(\mathbf{A}, \mathbf{B}) :=  \sqrt{\text{Dist}(\mathbf{A}, \mathbf{B}) \text{Dist}(\mathbf{B}, \mathbf{A})}$.
Nevertheless, none of these functions determine a metric on $M$, since the triangle inequality does not hold.}

Moreover, the dissipation inequality \eqref{decayCld}, which holds for
all relaxation processes, can be interpreted as follows:
during relaxation, the distance (measured in terms of energy)
between any solution $\mathbf{C}^{(2)}_{\text{i}}$ and a
constant solution $\mathbf{C}^{(1)}_{\text{i}} \equiv \mathbf{1}$ is not increasing.

\subsection{Sufficient condition for exponential stability}

Let us consider a loading program (strain-driven process)
$\{ \mathbf C \}_{t \in [0,T]}$ on the time interval $[0,T]$.
Let $\mathbf{C}^{(1)}_{\text{i}}$, $\mathbf{C}^{(2)}_{\text{i}} \in M$ be two solutions.
In order to prove the exponential stability, it is sufficient to prove that
there exists $t' \geq 0$ and $\gamma >0$
such that for all $t \geq t'$ (cf. \eqref{examp7})
\begin{equation}\label{assump1}
\big( \text{Dist}(\mathbf{C}^{(1)}_{\text{i}}, \mathbf{C}^{(2)}_{\text{i}})^2 \big)^{\cdot} \leq -
\gamma \ \text{Dist}(\mathbf{C}^{(1)}_{\text{i}}, \mathbf{C}^{(2)}_{\text{i}})^2.
\end{equation}
Indeed, in that case, using the Gronwall's inequality we get from \eqref{assump1} the following decay estimation
\begin{equation}\label{assump2}
\text{Dist}(\mathbf{C}^{(1)}_{\text{i}}(t), \mathbf{C}^{(2)}_{\text{i}}(t))^2 \leq
\text{Dist}(\mathbf{C}^{(1)}_{\text{i}}(t'), \mathbf{C}^{(2)}_{\text{i}}(t'))^2 \
\e^{-\gamma (t - t')}.
\end{equation}
Combining this result with \eqref{EnerDist3}, we get
the required estimation of type \eqref{decay}. Thus, the uniform
error estimation of Theorem 1 follows immediately from \eqref{assump2}.

\subsection{Reduction of the stability analysis to a
simplified problem with $\mathbf{C} = \mathbf{1}$}

Let $t^0$ be an arbitrary time instance.
In this section we discuss a procedure, which helps to simplify
the examination of the inequality \eqref{assump1} at time $t^0$.

The first
simplification of the problem is as follows.
We note that quantities \\ $\text{Dist}(\mathbf{C}^{(1)}_{\text{i}}(t^0),
\mathbf{C}^{(2)}_{\text{i}}(t^0))^2$, and
$\big( \text{Dist}(\mathbf{C}^{(1)}_{\text{i}}(t^0),
\mathbf{C}^{(2)}_{\text{i}}(t^0))^2 \big)^{\cdot}$
depend solely on $\overline{\mathbf{C}}(t^0)$, $\mathbf{C}^{(1)}_{\text{i}}(t^0)$, and
$\mathbf{C}^{(2)}_{\text{i}}(t^0)$ but not
on $\dot{\mathbf{C}}(t^0)$. Therefore, at the examination of \eqref{assump1} at $t=t^0$
 we can replace the actual loading programm
$\{ \mathbf C \}_{t \in [0,T]}$ by a constant loading (relaxation process):
we take a constant $\overline{\mathbf C} (t^0)$
instead of loading $\mathbf C(t)$,
where $ \overline{(\cdot)}$ stands for a unimodular part of a tensor.

The second simplification is as follows. Let
$\mathbf F_{0}$ be some "new" reference
configuration and $\det(\mathbf F_{0})=1 $. There is a one to one correspondence
between the solutions $\mathbf{C}^{(1)}_{\text{i}}(t)$, $\mathbf{C}^{(2)}_{\text{i}}(t)$
of the problem with the forcing function $\mathbf{C}(t)$ to
the solutions $(\mathbf{C}^{(1)}_{\text{i}})^{\text{new}}(t)$, $(\mathbf{C}^{(2)}_{\text{i}})^{\text{new}}(t)$ with
the forcing function $\mathbf{C}^{\text{new}}(t)$ (cf. Section 3.2)
\begin{equation}\label{simplific1}
\mathbf{C}^{\text{new}}(t) = \mathbf F^{-\text{T}}_{0} \mathbf{C} \mathbf F^{-1}_{0},
\quad
(\mathbf{C}^{(\text{k})}_{\text{i}})^{\text{new}}(t) =
\mathbf F^{-\text{T}}_{0} \mathbf{C}^{(\text{k})}_{\text{i}}(t) \mathbf F^{-1}_{0}, \quad
\text{k} \in \{1,2 \}.
\end{equation}
It follows from \eqref{EnerDist2} that
\begin{equation}\label{simplific2}
\text{Dist}\big( (\mathbf{C}^{(1)}_{\text{i}})^{\text{new}}(t^0),
(\mathbf{C}^{(2)}_{\text{i}})^{\text{new}}(t^0) \big)
 = \text{Dist} \big(\mathbf{C}^{(1)}_{\text{i}}(t^0), \mathbf{C}^{(2)}_{\text{i}}(t^0) \big),
\end{equation}
\begin{equation}\label{simplific3}
\big[ \text{Dist}\big( (\mathbf{C}^{(1)}_{\text{i}})^{\text{new}}(t),
(\mathbf{C}^{(2)}_{\text{i}})^{\text{new}}(t) \big)^2 \big]^{\cdot} |_{t=t^0}
 = \big[ \text{Dist} \big(\mathbf{C}^{(1)}_{\text{i}}(t),
 \mathbf{C}^{(2)}_{\text{i}}(t) \big)^2 \big]^{\cdot} |_{t=t^0},
\end{equation}
Therefore, estimation \eqref{assump1} is equivalent to
\begin{equation}\label{simplific4}
\big[ \text{Dist}\big( (\mathbf{C}^{(1)}_{\text{i}})^{\text{new}}(t),
(\mathbf{C}^{(2)}_{\text{i}})^{\text{new}}(t) \big)^2 \big]^{\cdot} |_{t=t^0} \leq -\gamma \
\text{Dist}\big( (\mathbf{C}^{(1)}_{\text{i}})^{\text{new}}(t^0),
(\mathbf{C}^{(2)}_{\text{i}})^{\text{new}}(t^0) \big)^2.
\end{equation}
Without loss of generality we assume $\det(\mathbf{C}(t^0))=1$.
By choosing $\mathbf F_{0} = \big(\mathbf C (t^0)\big)^{1/2}$
the problem can be reduced to the simplified problem with
$\mathbf C(t^0) = \mathbf{1}$.\footnote{Alternatively,
the problem can be reduced to the case $\mathbf{C}^{(1)}_{\text{i}}(t^0) = \mathbf{1}$
by choosing $\mathbf F_{0} = \big( \mathbf{C}^{(1)}_{\text{i}}(t^0)\big)^{1/2}$.}


\subsection{Evaluation of $\big( \text{Dist}(\mathbf{C}^{(1)}_{\text{i}}, \mathbf{C}^{(2)}_{\text{i}})^2 \big)^{\cdot}$}
In this section we evaluate $\big( \text{Dist}(\mathbf{C}^{(1)}_{\text{i}}, \mathbf{C}^{(2)}_{\text{i}})^2 \big)^{\cdot}$
at some fixed time instance $t^0$.
Without loss of generality (cf. the previous section)
it can be assumed that $\mathbf C(t^0) = \mathbf{1}$.
In that reduced case, the evolution equation \eqref{prob1} takes the form
\begin{equation}\label{evaluat0}
\dot{\mathbf{C}}_{\text{i}} =
\alpha( \| (\mathbf{C}^{-1}_{\text{i}})^{\text{D}} \| )
(\mathbf{C}^{-1}_{\text{i}})^{\text{D}} \ \mathbf{C}_{\text{i}}, \quad
\alpha(x) := \frac{\displaystyle 1}{\displaystyle
\eta \mu x}\Big\langle \frac{\displaystyle \mu x- \sqrt{2/3} K}{\displaystyle k_0}
\Big\rangle^{m}.
\end{equation}
Next, using the product rule we get from $(\ref{evaluat0})_1$
\begin{multline}\label{evaluat01}
\big(\mathbf{C}^{(1)}_{\text{i}} \mathbf{C}^{(2)-1}_{\text{i}} \big)^{\cdot} =
\dot{\mathbf{C}}^{(1)}_{\text{i}} \mathbf{C}^{(2)-1}_{\text{i}} +
\mathbf{C}^{(1)}_{\text{i}} (\mathbf{C}^{(2)-1}_{\text{i}})^{\cdot} \\ =
\mathbf{C}^{(1)}_{\text{i}} \mathbf{C}^{(1)-1}_{\text{i}} \dot{\mathbf{C}}^{(1)}_{\text{i}} \mathbf{C}^{(2)-1}_{\text{i}} +
\mathbf{C}^{(1)}_{\text{i}} (\mathbf{C}^{(2)-1}_{\text{i}})^{\cdot} \mathbf{C}^{(2)}_{\text{i}} \mathbf{C}^{(2)-1}_{\text{i}} \\
\stackrel{(\ref{evaluat0})_1}{=} \mathbf{C}^{(1)}_{\text{i}} \big[ \alpha^{(1)} (\mathbf{C}^{(1)-1}_{\text{i}})^{\text{D}} -
\alpha^{(2)} (\mathbf{C}^{(2)-1}_{\text{i}})^{\text{D}} \big] \mathbf{C}^{(2)-1}_{\text{i}},
\end{multline}
where $\alpha^{(k)} := \alpha (\| (\mathbf{C}^{(k)-1}_{\text{i}})^{\text{D}} \|)$ for $k \in \{1,2\}$.

Further, we compute the derivative of
$\psi_{\text{el}} (\mathbf{A})$ using a coordinate-free tensor setting (see, for example, \cite{Itskov, Shutov2}).
\begin{equation}\label{evaluat}
\rho_{\scriptscriptstyle \text{R}} \frac{\displaystyle \partial \psi_{\text{el}} (\mathbf{A})}{\displaystyle \partial \mathbf{A}} =
\frac{k}{2} \ln \sqrt{\det \mathbf{A}}
\ \mathbf{A}^{- \text{T}} + \frac{\mu}{2}
\mathbf{A}^{- \text{T}} \ (\overline{\mathbf{A}}^{\text{T}})^{\text{D}}.
\end{equation}
We abbreviate $\Delta := \| \mathbf{C}^{(2)}_{\text{i}} - \mathbf{C}^{(1)}_{\text{i}} \|$.
Note that
(see Appendix B), since $\mathbf{C}^{(1)}_{\text{i}}, \mathbf{C}^{(2)}_{\text{i}} \in M$
\begin{equation}\label{evaluat30}
\text{tr} \big( ( \mathbf{C}^{(2)-1}_{\text{i}} -
\mathbf{C}^{(1)-1}_{\text{i}}) \mathbf{C}^{(1)}_{\text{i}}
\big) = O (\Delta^2), \quad
\text{tr} \big( \mathbf{C}^{(1)}_{\text{i}}
( \mathbf{C}^{(1)-1}_{\text{i}} - \mathbf{C}^{(2)-1}_{\text{i}})
\big) = O (\Delta^2),
\text{as} \ \Delta \rightarrow 0.
\end{equation}
Thus, using $(\ref{evaluat30})_1$ we get
\begin{equation}\label{evaluat3}
\big( \mathbf{C}^{(2)-1}_{\text{i}} \mathbf{C}^{(1)}_{\text{i}}\big )^{\text{D}} =
\big( (\mathbf{C}^{(2)-1}_{\text{i}} - \mathbf{C}^{(1)-1}_{\text{i}})  \mathbf{C}^{(1)}_{\text{i}} + \mathbf{1} \big)^{\text{D}} =
(\mathbf{C}^{(2)-1}_{\text{i}} - \mathbf{C}^{(1)-1}_{\text{i}})  \mathbf{C}^{(1)}_{\text{i}} + O (\Delta^2).
\end{equation}
Combining \eqref{evaluat} with \eqref{evaluat3} we get
\begin{multline}\label{evaluat4}
\rho_{\scriptscriptstyle \text{R}} \frac{\displaystyle \partial \psi_{\text{el}} (\mathbf{C}^{(1)}_{\text{i}} \mathbf{C}^{(2)-1}_{\text{i}})}
{\displaystyle \partial (\mathbf{C}^{(1)}_{\text{i}} \mathbf{C}^{(2)-1}_{\text{i}})} =
\frac{\mu}{2}
\big(\mathbf{C}^{(1)}_{\text{i}} \mathbf{C}^{(2)-1}_{\text{i}}\big)^{- \text{T}} \
\big[\overline{\big(\mathbf{C}^{(1)}_{\text{i}} \mathbf{C}^{(2)-1}_{\text{i}}\big)}^{\text{T}} \big]^{\text{D}} \\ =
\frac{\mu}{2} \mathbf{C}^{(1)-1}_{\text{i}} \mathbf{C}^{(2)}_{\text{i}}
(\mathbf{C}^{(2)-1}_{\text{i}} - \mathbf{C}^{(1)-1}_{\text{i}})  \mathbf{C}^{(1)}_{\text{i}} + O (\Delta^2)
= \frac{\mu}{2}
(\mathbf{C}^{(2)-1}_{\text{i}} - \mathbf{C}^{(1)-1}_{\text{i}})  \mathbf{C}^{(1)}_{\text{i}} + O (\Delta^2).
\end{multline}
Next, denote by $\acute{\alpha}$ the derivative of $\alpha(x)$
 at $x=\| (\mathbf{C}^{(1)-1}_{\text{i}})^{\text{D}} \|$.
Therefore,
\begin{multline}\label{alphaDeriv}
\alpha^{(1)} - \alpha^{(2)} = \acute{\alpha} \
(\| (\mathbf{C}^{(1)-1}_{\text{i}})^{\text{D}} \| -
\| (\mathbf{C}^{(2)-1}_{\text{i}})^{\text{D}} \| ) + O (\Delta^2) \\ =
\frac{\acute{\alpha}}{\| (\mathbf{C}^{(1)-1}_{\text{i}})^{\text{D}} \|}
(\mathbf{C}^{(1)-1}_{\text{i}})^{\text{D}} : (\mathbf{C}^{(1)-1}_{\text{i}} - \mathbf{C}^{(2)-1}_{\text{i}})
+ O (\Delta^2).
\end{multline}

It can be assumed that the
overstress $f = \mu \| (\mathbf{C}^{(1)-1}_{\text{i}})^{\text{D}} \| - \sqrt{\frac{2}{3}} K$
is bounded by $\sqrt{\frac{2}{3}} K$.
Thus, we suppose
$ \sqrt{\frac{2}{3}} K/\mu  <  \| (\mathbf{C}^{(1)-1}_{\text{i}})^{\text{D}} \| \leq 2 \sqrt{\frac{2}{3}} K/\mu $.
Here, the first inequality is needed to ensure the overstress
is larger than zero.

Using the property $\mathbf{A}:(\mathbf{B} \mathbf{C} \mathbf{D}) =
(\mathbf{B}^{\text{T}} \mathbf{A} \mathbf{D}^{\text{T}}) : \mathbf{C}$ it
follows from \eqref{evaluat01} and \eqref{evaluat4} that
\begin{multline}\label{evaluat5}
\big( \text{Dist}(\mathbf{C}^{(1)}_{\text{i}}, \mathbf{C}^{(2)}_{\text{i}})^2 \big)^{\cdot} \stackrel{\eqref{EnerDist1}}{=}
\big( \rho_{\scriptscriptstyle \text{R}} \psi_{\text{el}} (\mathbf{C}^{(1)}_{\text{i}} \mathbf{C}^{(2)-1}_{\text{i}}) \big)^{\cdot}
= \rho_{\scriptscriptstyle \text{R}} \frac{\displaystyle \partial \psi_{\text{el}} (\mathbf{C}^{(1)}_{\text{i}} \mathbf{C}^{(2)-1}_{\text{i}})}
{\displaystyle \partial (\mathbf{C}^{(1)}_{\text{i}} \mathbf{C}^{(2)-1}_{\text{i}})} :
\big(\mathbf{C}^{(1)}_{\text{i}} \mathbf{C}^{(2)-1}_{\text{i}} \big)^{\cdot}
\\ \stackrel{\eqref{evaluat01}, \eqref{evaluat4}}{=}
-\frac{\mu}{2} \big(\mathbf{C}^{(1)}_{\text{i}}
(\mathbf{C}^{(1)-1}_{\text{i}} - \mathbf{C}^{(2)-1}_{\text{i}})  \big):\big( \alpha^{(1)} (\mathbf{C}^{(1)-1}_{\text{i}})^{\text{D}} -
\alpha^{(2)} (\mathbf{C}^{(2)-1}_{\text{i}})^{\text{D}} \big) + O(\Delta^3) \\
\stackrel{(\ref{evaluat30})_2}{=}
-\frac{\mu}{2} \big(\mathbf{C}^{(1)}_{\text{i}}
(\mathbf{C}^{(1)-1}_{\text{i}} - \mathbf{C}^{(2)-1}_{\text{i}})  \big):\big( \alpha^{(1)} \mathbf{C}^{(1)-1}_{\text{i}} -
\alpha^{(2)} \mathbf{C}^{(2)-1}_{\text{i}} \big) + O(\Delta^3) \\
= -\frac{\mu}{2} \alpha^{(1)} \big(\mathbf{C}^{(1)}_{\text{i}}
(\mathbf{C}^{(1)-1}_{\text{i}} - \mathbf{C}^{(2)-1}_{\text{i}})  \big):
\big( \mathbf{C}^{(1)-1}_{\text{i}} -
\mathbf{C}^{(2)-1}_{\text{i}} \big)   \\
 -\frac{\mu}{2} (\alpha^{(1)}- \alpha^{(2)}) \big(\mathbf{C}^{(1)}_{\text{i}}
(\mathbf{C}^{(1)-1}_{\text{i}} - \mathbf{C}^{(2)-1}_{\text{i}})  \big): \mathbf{C}^{(2)-1}_{\text{i}} + O(\Delta^3)
\stackrel{\eqref{alphaDeriv}}{=}
F_{I} + F_{II} + O(\Delta^3),
\end{multline}
where $F_{I}$ and $F_{II}$ are given by
\begin{equation}\label{evaluat6}
F_{I} := -\frac{\mu}{2} \alpha^{(1)} \text{tr} \Big( \big( \mathbf{C}^{(1)-1}_{\text{i}} -
\mathbf{C}^{(2)-1}_{\text{i}} \big)  \mathbf{C}^{(1)}_{\text{i}}   \big( \mathbf{C}^{(1)-1}_{\text{i}} -
\mathbf{C}^{(2)-1}_{\text{i}} \big) \Big),
\end{equation}
\begin{equation}\label{evaluat7}
F_{II} := -\frac{\mu}{2} \frac{\acute{\alpha}}{\| (\mathbf{C}^{(1)-1}_{\text{i}})^{\text{D}} \|} \big(
(\mathbf{C}^{(1)-1}_{\text{i}})^{\text{D}} : (\mathbf{C}^{(1)-1}_{\text{i}} - \mathbf{C}^{(2)-1}_{\text{i}}) \big)
\big( \mathbf{1} : (\mathbf{C}^{(1)-1}_{\text{i}} - \mathbf{C}^{(2)-1}_{\text{i}}) \big).
\end{equation}

Now, for any pair of real positive numbers $(\theta, \Delta)$ let us define a subset of $M \times M$ by
\begin{equation}\label{evaluat8}
S(\theta, \Delta) :=
\{ (\mathbf{C}^{(1)}_{\text{i}}, \mathbf{C}^{(2)}_{\text{i}}) \in M \times M \ | \
 \| (\mathbf{C}^{(1)-1}_{\text{i}})^{\text{D}} \| \leq \theta, \
 \| \mathbf{C}^{(1)}_{\text{i}} -  \mathbf{C}^{(2)}_{\text{i}} \| \leq \Delta  \}.
\end{equation}

By definition, put
\begin{equation}\label{evaluat9}
\Phi(\mathbf{C}^{(1)}_{\text{i}}, \mathbf{C}^{(2)}_{\text{i}} ) := - \frac{2 \alpha^{(1)} F_{II}}{\acute{\alpha} F_{I}} =
\frac{-2 \big( \frac{(\mathbf{C}^{(1)-1}_{\text{i}})^{\text{D}}}{\| (\mathbf{C}^{(1)-1}_{\text{i}})^{\text{D}} \|}
 : (\mathbf{C}^{(1)-1}_{\text{i}} - \mathbf{C}^{(2)-1}_{\text{i}}) \big)
\big( \mathbf{1} : (\mathbf{C}^{(1)-1}_{\text{i}} - \mathbf{C}^{(2)-1}_{\text{i}}) \big)}{\text{tr}
\Big( \big( \mathbf{C}^{(1)-1}_{\text{i}} -
\mathbf{C}^{(2)-1}_{\text{i}} \big)  \mathbf{C}^{(1)}_{\text{i}}   \big( \mathbf{C}^{(1)-1}_{\text{i}} -
\mathbf{C}^{(2)-1}_{\text{i}} \big) \Big)}.
\end{equation}


There exists a function  $q(\theta)>0$ such that
\begin{equation}\label{evaluat10}
q(\theta) \geq \Phi(\mathbf{C}^{(1)}_{\text{i}}, \mathbf{C}^{(2)}_{\text{i}} ) + O(\Delta) \quad
\text{for all} \ (\mathbf{C}^{(1)}_{\text{i}}, \mathbf{C}^{(2)}_{\text{i}} )
\in  S(\theta, \Delta).
\end{equation}

The numerical evaluation of the function $q(\theta)$ is discussed in the Appendix C.
Moreover, suppose that
\begin{equation}\label{evaluat11}
\alpha^{(1)} \geq q \Big( 2 \sqrt{\frac{2}{3}} K/\mu  \Big) \acute{\alpha}.
\end{equation}
This condition will be discussed in the next section.
Multiplying both sides of \eqref{evaluat11} by $F_{I} \frac{\displaystyle 1}{\displaystyle \alpha^{(1)}} < 0$
and noting that $\frac{\displaystyle \acute{\alpha}}{\displaystyle \alpha^{(1)}} F_{I} O(\Delta) = O(\Delta^3)$ we get
for all $(\mathbf{C}^{(1)}_{\text{i}}, \mathbf{C}^{(2)}_{\text{i}} )
\in  S(2 \sqrt{\frac{2}{3}} K/\mu, \Delta)$
\begin{equation}\label{evaluat122}
F_{I} \stackrel{\eqref{evaluat11}}{\leq} q \Big(2 \sqrt{\frac{2}{3}} K/\mu \Big) \frac{\displaystyle \acute{\alpha}}{\displaystyle \alpha^{(1)}} F_{I}
\stackrel{\eqref{evaluat10}}{\leq}
\big(\Phi(\mathbf{C}^{(1)}_{\text{i}}, \mathbf{C}^{(2)}_{\text{i}} ) + O(\Delta) \big)
\frac{\displaystyle \acute{\alpha}}{\displaystyle \alpha^{(1)}} F_{I}
 \stackrel{\eqref{evaluat9}}{=} -2 F_{II} + O(\Delta^3).
\end{equation}
Multiplying both sides of \eqref{evaluat122} by $1/2$ and adding
$1/2 F_{I} + F_{II}$, we get 
\begin{equation}\label{evaluat123}
F_{I}+F_{II} \leq 1/2 F_{I} + O(\Delta^3).
\end{equation}
Combining this result with \eqref{evaluat5} we obtain
\begin{equation}\label{evaluat13}
\big( \text{Dist}(\mathbf{C}^{(1)}_{\text{i}},
\mathbf{C}^{(2)}_{\text{i}})^2 \big)^{\cdot} \leq 1/2 F_{I} + O(\Delta^3)
.
\end{equation}
Next, if $f \geq f_0$ for some $f_0 >0$, then there exists $C_5>0$ such that
\begin{equation}\label{evaluat14}
F_{I} = -\frac{\mu}{2} \alpha^{(1)}
\big \| ( \mathbf{C}^{(1)-1}_{\text{i}} -  \mathbf{C}^{(2)-1}_{\text{i}} )
(\mathbf{C}^{(1)}_{\text{i}})^{1/2} \big \|^2 \leq - C_5 \Delta^2.
\end{equation}
Therefore, for small $\Delta$,
inequality \eqref{evaluat13} yields
\begin{equation}\label{evaluat1334}
\big( \text{Dist}(\mathbf{C}^{(1)}_{\text{i}}, \mathbf{C}^{(2)}_{\text{i}})^2 \big)^{\cdot} \leq 1/4 F_{I}.
\end{equation}
Similarly to the proof of \eqref{EnerDist3} we obtain with some $C_6 > 0$
\begin{multline}\label{evaluat142}
F_{I} = -\frac{\mu}{2} \alpha^{(1)}
\big \| ( \mathbf{C}^{(1)-1}_{\text{i}} -  \mathbf{C}^{(2)-1}_{\text{i}} )
(\mathbf{C}^{(1)}_{\text{i}})^{1/2} \big \|^2 \\
\leq -\frac{\mu}{2} \alpha^{(1)} C_6 \big \| (\mathbf{C}^{(1)}_{\text{i}})^{1/2} \big \|^{-2} \
 \big \| (\mathbf{C}^{(1)}_{\text{i}})^{1/2} ( \mathbf{C}^{(1)-1}_{\text{i}} -  \mathbf{C}^{(2)-1}_{\text{i}} )
(\mathbf{C}^{(1)}_{\text{i}})^{1/2} \big \|^2 \\
=
-\frac{\mu}{2} \alpha^{(1)} C_6 \big \| (\mathbf{C}^{(1)}_{\text{i}})^{1/2} \big \|^{-2} \
 \big \| (\mathbf{C}^{(1)}_{\text{i}})^{1/2} \mathbf{C}^{(2)-1}_{\text{i}}
 (\mathbf{C}^{(1)}_{\text{i}})^{1/2}   - \mathbf{1} \big \|^2 \\
 \leq -\frac{\mu}{2} \alpha^{(1)} ( C_6 / C'_4 ) \big \| (\mathbf{C}^{(1)}_{\text{i}})^{1/2} \big \|^{-2}
 \rho_{\scriptscriptstyle \text{R}} \psi_{\text{el}}
 ((\mathbf{C}^{(1)}_{\text{i}})^{1/2} \mathbf{C}^{(2)-1}_{\text{i}}
 (\mathbf{C}^{(1)}_{\text{i}})^{1/2}) \\
= -\frac{\mu}{2} \alpha^{(1)} ( C_6 / C'_4 ) \big \| (\mathbf{C}^{(1)}_{\text{i}})^{1/2} \big \|^{-2}
\rho_{\scriptscriptstyle \text{R}} \psi_{\text{el}}
(\mathbf{C}^{(1)}_{\text{i}} \mathbf{C}^{(2)-1}_{\text{i}}).
\end{multline}
Finally, combining \eqref{evaluat1334} with
\eqref{evaluat142} we get the required estimation \eqref{assump1} if the following assumptions hold:
$ 0 < f_0 \leq f \leq \sqrt{\frac{2}{3}} K$, $\alpha^{(1)} \geq q \Big( 2 \sqrt{\frac{2}{3}} K/\mu  \Big) \acute{\alpha}$.

\subsection{Analysis of the sufficient stability condition}
In this section we analyze the condition \eqref{evaluat11} which was used in the previous section
to prove the inequality \eqref{assump1}.
First, we suppose $\| (\mathbf{C}^{-1}_{\text{i}})^{\text{D}} \| > \sqrt{\frac{2}{3}} K/\mu  $
to ensure the overstress is larger than zero.
Using $(\ref{evaluat0})_2$ it can be easily shown that
\eqref{evaluat11} is equivalent to
\begin{equation}\label{evaluat16}
\| (\mathbf{C}^{-1}_{\text{i}})^{\text{D}} \| \geq x_{cr},
\end{equation}
where the critical value $x_{cr}$ is given by
\begin{equation}\label{evaluat170}
x_{cr} := \Big( \sqrt{\frac{2}{3}} K + (m-1) \mu q(\theta) +
\sqrt{ (\sqrt{\frac{2}{3}} K + (m-1) \mu q(\theta))^2 + 4 \mu q(\theta) \sqrt{\frac{2}{3}} K }
\Big)/\big(2 \mu \big),
\end{equation}
with $\theta = 2 \sqrt{\frac{2}{3}} K/\mu$.
For small values of $q(2 \sqrt{\frac{2}{3}} K/\mu)$, a simple estimation for $x_{cr}$  is valid
\begin{equation}\label{evaluat18}
x_{cr} = \sqrt{\frac{2}{3}} K / \mu + m q(2 \sqrt{\frac{2}{3}} K/\mu) + O \big( (q(2 \sqrt{\frac{2}{3}} K/\mu))^2\big).
\end{equation}
Alternatively, in terms of the overstress $f$, the condition \eqref{evaluat11} is equivalent to
\begin{equation}\label{evaluat17}
f \geq f_{cr},
\end{equation}
where the critical overstress $f_{cr}$ is estimated by
\begin{equation}\label{evaluat182}
f_{cr} =  m \mu q(2 \sqrt{\frac{2}{3}} K/\mu) + O \big( (q(2 \sqrt{\frac{2}{3}} K/\mu))^2\big).
\end{equation}

\begin{figure}\centering
\psfrag{A}[m][][1][0]{$-\sqrt{\frac{2}{3}} K$}
\psfrag{B}[m][][1][0]{$0$}
\psfrag{C}[m][][1][0]{$f_{cr}$}
\psfrag{D}[m][][1][0]{$\sqrt{\frac{2}{3}} K$}
\psfrag{E}[m][][1][0]{$f$}
\psfrag{F}[m][][1][0]{No flow}
\psfrag{G}[m][][1][0]{Exponential stability}
\scalebox{0.9}{\includegraphics{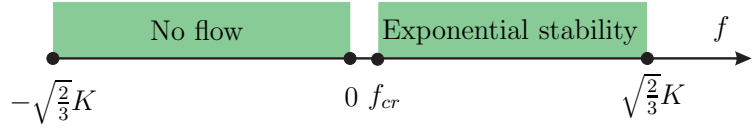}}
\caption{Domain of exponential stability. \label{fig3}}
\end{figure}

The situation is summarized in figure \ref{fig3}.

For instance, for aluminium alloy we put $K=300$ MPa, $\mu = 25 000$ MPa.
Thus $ 2 \sqrt{\frac{2}{3}} K/ \mu \approx 0.014$.
Next, $q(0.014) \approx 0.00000023$ (See Appendix C).
 Therefore, the critical overstress is given by
$f_{cr} \approx  m \ 0.0057$ MPa. For physically reasonable values of $m$ ($m \leq 100$)
this critical value is negligible compared to the size of the elastic domain
$\sqrt{\frac{2}{3}} K \approx 245$ MPa.

\textbf{Remark.}
Since the overstress $f$ is isolated from zero due to the sufficient stability
condition \eqref{evaluat17}, the current theory can not be applied to
exactly quasistatic processes. On the other hand, the
theory is directly applicable for \emph{nearly quasistatic} processes with the oversress $f$ larger that
$f_{cr}$.

\section{Accuracy testing of implicit integrators}

The numerical implementation of the material model \eqref{prob1} ---  \eqref{prob3} within a
displacement based Finite Element Method (FEM) with implicit time stepping
is based on the implicit integration of the evolution equation \eqref{prob1}
(see, for example, \cite{SimHug}).
This procedure should provide the stresses
as a function of the strain history.

More precisely, suppose that the right Cauchy-Green tensor ${}^{n+1} \mathbf C$
at the time $t_{n+1} = t_n + \Delta t$
is known and assume that the internal variable
$\mathbf C_{\text{i}}$ at the time $t_n$ is given by
${}^{n} \mathbf C_{\text{i}}$.
We need to compute the internal variable $\mathbf C_{\text{i}}$ at the time $t_{n+1}$ in order
to evaluate the stress tensor
${}^{n+1} \tilde{\mathbf T} =
\tilde{\mathbf T} ({}^{n+1} \mathbf C, {}^{n+1}
\mathbf C_{\text{i}})$.

Note that the norm of the driving force $\mathfrak{F}$ and the overstress $f$
can be represented as functions of ${}^{n+1} {\mathbf C}$ and ${}^{n+1} \mathbf C_{\text{i}}$:
\begin{equation}\label{IntrodOfB2}
\mathfrak{F}({}^{n+1} {\mathbf C}, {}^{n+1} \mathbf C_{\text{i}})= \sqrt{\text{tr}
\big[ \big( {}^{n+1} \mathbf C \ \tilde{\mathbf T} ({}^{n+1}
{\mathbf C}, {}^{n+1} \mathbf C_{\text{i}}) \big)^{\text{D}} \big]^2 },
\end{equation}
\begin{equation}\label{Per0}
f({}^{n+1} {\mathbf C}, {}^{n+1} \mathbf C_{\text{i}}) =
\mathfrak{F}({}^{n+1} {\mathbf C}, {}^{n+1} \mathbf C_{\text{i}}) - \sqrt{\frac{2}{3}}K.
\end{equation}

For what follows it is useful to introduce the incremental inelastic parameter
\begin{equation}\label{ineinc}
\xi := \Delta t \ {}^{n+1} \lambda_{\text{i}}.
\end{equation}
Thus, according to the Perzyna rule, we get the following equation with respect to
${}^{n+1} {\mathbf C}$, ${}^{n+1} \mathbf C_{\text{i}}$ and $\xi$
\begin{equation}\label{Per2}
\xi= \frac{\displaystyle \Delta t}{\displaystyle \eta}\Big \langle
\frac{\displaystyle f({}^{n+1} {\mathbf C}, {}^{n+1} \mathbf C_{\text{i}})}{\displaystyle k_0} \Big\rangle^m.
\end{equation}
The remaining equation for finding unknown ${}^{n+1} \mathbf C_{\text{i}}$ and $\xi$ is obtained through
the time discretization of \eqref{prob1}, which will be discussed in the next section.

\subsection{Euler Backward Method and geometric implicit integrators}

We introduce a nonlinear operator
$\mathbf B ({}^{n+1} {\mathbf C}, {}^{n+1} \mathbf C_{\text{i}}, \xi)$ as
\begin{equation}\label{IntrodOfB}
\mathbf B ({}^{n+1} {\mathbf C}, {}^{n+1} {\mathbf C}_{\text{i}}, \xi) :=
2 \frac{\xi}{\mathfrak{F} ({}^{n+1} {\mathbf C}, {}^{n+1} \mathbf C_{\text{i}})}
\big({}^{n+1} \mathbf C \ \tilde{\mathbf T} ({}^{n+1}
{\mathbf C}, {}^{n+1} \mathbf C_{\text{i}}) \big)^{\text{D}}.
\end{equation}

Let us consider the classical Euler-Backward method (EBM)
(see, for example, \cite{DettRes, Hartmann, SimHug})
 being applied to the evolution problem \eqref{prob1}

\begin{equation}\label{classicalEuler}
{}^{n+1} {\mathbf C}_{\text{i}} =
 \big[ \mathbf 1 - \mathbf B ({}^{n+1} {\mathbf C}, {}^{n+1}
 {\mathbf C}_{\text{i}}, \xi) \big]^{-1} \
    {}^{n} {\mathbf C}_{\text{i}}.
\end{equation}
Since the symmetry of the internal variable ${}^{n+1} {\mathbf C}_{\text{i}}$
is exactly preserved by the EBM\footnote{
Moreover, it was shown in \cite{Shutov1}
that the \emph{symmetry is exactly preserved} by Euler-Backward method and Exponential Method
even in a more general case of a nonlinear kinematic hardening.},
this equation is equivalent to
\begin{equation}\label{classicalEuler2}
{}^{n+1} {\mathbf C}_{\text{i}} = \text{sym} \big(
 \big[ \mathbf 1 - \mathbf B ({}^{n+1} {\mathbf C}, {}^{n+1} {\mathbf C}_{\text{i}}, \xi) \big]^{-1} \
    {}^{n} {\mathbf C}_{\text{i}}\big).
\end{equation}
The modified Euler-Backward method (MEBM) (see \cite{Helm2, Shutov1}) uses the following equation
\begin{equation}\label{MEBM1}
{}^{n+1} {\mathbf C}_{\text{i}} = \overline{\text{sym} \big(
 \big[ \mathbf 1 - \mathbf B ({}^{n+1} {\mathbf C}, {}^{n+1} {\mathbf C}_{\text{i}}, \xi) \big]^{-1} \
    {}^{n} {\mathbf C}_{\text{i}}\big)}.
\end{equation}
Finally, the Exponential Method (EM) (see, for instance,
\cite{DettRes, MiStei, Miehe, WebAnan} )
is based on the use of the tensor exponential $\exp (\cdot)$.
As it was shown in \cite{Shutov1}, the Exponential Method
can be written in the following form:
\begin{equation}\label{EM1}
{}^{n+1} {\mathbf C}_{\text{i}} = \overline{\text{sym} \big(
 \exp \big[  \mathbf B ({}^{n+1} {\mathbf C}, {}^{n+1}
 {\mathbf C}_{\text{i}}, \xi) \big] \
    {}^{n} {\mathbf C}_{\text{i}}\big)}.
\end{equation}

Combining \eqref{Per2} with one of
the discretization methods (equations \eqref{classicalEuler2}, \eqref{MEBM1} or \eqref{EM1})
a closed system of equations is obtained.
One possible solution strategy for the resulting problem
was discussed in \cite{Shutov1}, and the application of a coordinate-free
tensor formalism to the numerical solution was analyzed in \cite{Shutov2}.

We note that the geometric property of the exact flow
($\mathbf{C}_{\text{i}} \in M$) is
exactly satisfied by MEBM and EM. Therefore we refer to these two methods as to
\emph{geometric integrators}. On the other hand, the incompressibility constraint
is violated by the classical EBM.

For all the three methods, the error on the step is
bounded by the second power of the step size (cf. estimation \eqref{ErrStep}),
if the right-hand side is a smooth function.
Strong local nonlinearities due to the distinction
into elastic and inelastic material behavior or
due to the non-smoothness of the loading function $\mathbf{C} (t)$ may
increase the error on the step.


\subsection{Testing results}

The theoretical
results obtained in this study are validated
via a series of numerical tests.
Let us simulate the material behavior under strain controlled,
nonproportional and non-monotonic loading
in the time interval $t \in [0,300]$.
Suppose that the deformation gradient is defined by
\begin{equation}\label{loaprog0}
\mathbf F (t) = \overline{\mathbf F^{\prime} (t)},
\end{equation}
where $\mathbf F^{\prime}(t)$ is a piecewise linear function of time $t$ such that
$\mathbf F^{\prime} (0) = \mathbf F_1$, $\mathbf F^{\prime} (100) = \mathbf F_2$,
$\mathbf F^{\prime} (200) = \mathbf F_3$, and $\mathbf F^{\prime} (300) = \mathbf F_4$ with
\begin{equation*}\label{loaprog2}
\mathbf F_1 :=\mathbf 1, \
\mathbf F_2 := \left(
\begin{array}{ccc}
2 &  & 0 \\
0 & \displaystyle \frac{1}{\sqrt2} & 0 \\
0 & 0 & \displaystyle \frac{1}{\sqrt2}
\end{array}
\right), \
\mathbf F_3 := \left(
\begin{array}{ccc}
1 & 1  & 0 \\
0 & 1 & 0 \\
0 & 0 & 1
\end{array}
\right), \
\mathbf F_4 := \left(
\begin{array}{ccc}
\displaystyle \frac{1}{\sqrt2} &  & 0 \\
0 & 2 & 0 \\
0 & 0 & \displaystyle \frac{1}{\sqrt2}
\end{array}
\right).
\end{equation*}
Thus, we put
\begin{equation*}\label{loaprog}
\mathbf F^{\prime} (t) :=
\begin{cases}
    (1 - t/100) \mathbf F_1  + (t/100) \mathbf F_2 \quad \text{if} \ t \in [0,100] \\
    (2 - t/100) \mathbf F_2  + (t/100-1) \mathbf F_3 \quad \text{if} \ t \in (100,200] \\
     (3 - t/100) \mathbf F_3  + (t/100-2) \mathbf F_4 \quad \text{if} \ t \in (200,300]
\end{cases}.
\end{equation*}

The material parameters used in simulations are summarized in table \ref{table2}.

\begin{figure}\centering
\psfrag{A}[m][][1][0]{$t [s]$}
\psfrag{B}[m][][1][0]{$ \| {\mathbf C}^{numer}_{\text{i}} - {\mathbf C}^{exact}_{\text{i}} \|  $}
\psfrag{C}[m][][1][0]{MEBM and EM}
\psfrag{D}[m][][1][0]{EBM}
\psfrag{E}[m][][1][0]{$\Delta t = 1$ s}
\psfrag{F}[m][][1][0]{$\Delta t = 0.5$ s}
\scalebox{0.8}{\includegraphics{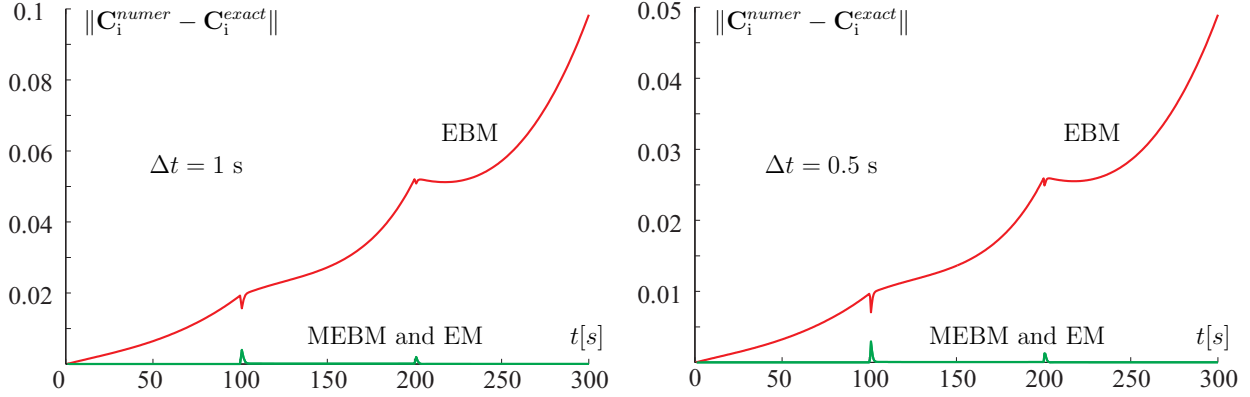}}
\caption{Accuracy analysis concerning
Euler Backward Method (EBM),
Modified Euler Backward Method (MEBM),
and Exponential Method (EM). \label{fig4}}
\end{figure}

\begin{table}[h]
\caption{Material parameters}
\begin{tabular}{| l l l l l l |}
\hline
$k$ [MPa] & $\mu$ [MPa]  & $K$ [MPa] & $m$ [-] & $\eta$ [$\text{s}^{-1}$] & $k_0$ [Mpa]  \\ \hline
73500 & 28200  &   270 & 3.6 &  $2 \cdot 10^6$ & 1 \\ \hline
\end{tabular} 
\label{table2}
\end{table}
Next, we suppose that the reference configuration is stress free.
 Therefore we put
\begin{equation}\label{Inico}
\mathbf C_{\text{i}}|_{t=0} =  \mathbf 1.
\end{equation}
The numerical solution obtained with extremely small time step
($\Delta t = 0.01 \text{s}$) will be named the \emph{exact solution} and denoted by
${\mathbf C}^{exact}_{\text{i}}$. Next, the numerical solutions with
$\Delta t = 1 \text{s}$ and $\Delta t = 0.5 \text{s}$ are
denoted by ${\mathbf C}^{numer}_{\text{i}}$. The error
$\| {\mathbf C}^{numer}_{\text{i}} - {\mathbf C}^{exact}_{\text{i}} \|$
is plotted on figure \ref{fig4}.

For all three methods the error is proportional to $\Delta t$.
Moreover, in accordance with Theorem 1 (cf. Section 2.2), the error is \emph{uniformly}
bounded for geometric integrators (MEBM and EM). More precisely,
the error is bounded by $C \Delta t$, where the constant $C$ \emph{does not depend} on the
size of the entire time interval.
Next, since the incompressibility condition is violated by EBM, the geometric property \eqref{geopr0}
is lost and some \emph{spurious degrees of freedom} are introduced. In that case,
only a weaker error estimation is valid: $\| {\mathbf C}^{numer}_{\text{i}} - {\mathbf C}^{exact}_{\text{i}} \| \leq
\widetilde{C}(T) \Delta t$, where $\widetilde{C}(T)$ depends on the size $T$ of the entire time interval.

\section{Discussion and conclusion}

\begin{figure}\centering
\scalebox{1.0}{\includegraphics{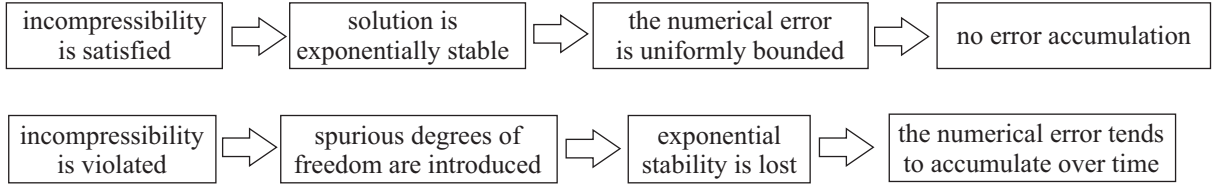}}
\caption{Summary of the main results. \label{fig6}}
\end{figure}

In the last decade, intensive research has been carried out
concerning the development of so-called geometric integrators for the
evolution equations of finite plasticity/viscoplasticity, which exactly preserve
the inelastic incompressibility condition. The excellent accuracy and convergence
properties of such algorithms were analyzed by numerical computations.
Particularly, the long term accuracy of geometric integrators
was analyzed in the paper \cite{Shutov1}, and the
\emph{absence of error accumulation} was numerically verified.
In the current study, a rigorous mathematical formulation
of this phenomena is proposed. The main result of the current paper is as follows:
the numerical error is \emph{uniformly} bounded by $C \Delta t$ if
the incompressibility condition is satisfied.
In terms of a classical model of finite viscoplasticity
we prove that \emph{all} first order accurate geometric integrators
are \emph{equivalent} in that sense. This theoretical result
corresponds with the numerical tests. Indeed, MEBM and EM are equivalent
concerning the accuracy and convergence (cf. figure \ref{fig4}).
The main results are summarized diagrammatically on figure \ref{fig6}.

The property of the exponential stability of the exact
plastic flow was mathematically analyzed in this paper.
Obviously, that property must be utilized
during the development of new material models and
corresponding algorithms in order to improve the accuracy
and convergence of numerical computations.

\section*{Appendix A}

Suppose $\Delta = \| \mathbf{\Delta} \| \rightarrow 0$.
Let us show that
\begin{equation}\label{EnerDist52}
\rho_{\scriptscriptstyle \text{R}} \psi_{\text{el}} (\mathbf{1} + \mathbf{\Delta}) =
\frac{k}{8} (\text{tr} \mathbf{\Delta})^2 +
\frac{\mu}{4} \text{tr} \big( (\mathbf{\Delta}^{\text{D}} )^2 \big) + O( \Delta^3 ).
\end{equation}

First, recall the Taylor expansion of $\det (\mathbf{1} + \mathbf{\Delta})$ up to second order
\begin{equation}\label{TaylorDet}
\det (\mathbf{1} + \mathbf{\Delta}) = 1 + \text{tr}(\mathbf{\Delta}) + 1/2 (\text{tr}(\mathbf{\Delta}))^2
-1/2\text{tr}(\mathbf{\Delta}^2) + O( \Delta^3 ).
\end{equation}
Therefore,
\begin{equation}\label{AppA2}
\sqrt{\det (\mathbf{1} + \mathbf{\Delta})} = 1 + 1/2\text{tr}(\mathbf{\Delta}) + O( \Delta^2),
\end{equation}
\begin{equation}\label{AppA3}
\frac{k}{2}\big(\text{ln}\sqrt{\det (\mathbf{1} + \mathbf{\Delta})} \big)^2
=  \frac{k}{8} (\text{tr} \mathbf{\Delta})^2 + O( \Delta^3).
\end{equation}
Next, note that for small $\varepsilon$ we have
\begin{equation}\label{AppA33}
(1 + \varepsilon)^{-1/3}= 1 - 1/3 \varepsilon + 2/9 \varepsilon^2 + O (\varepsilon^3).
\end{equation}
Combining this with \eqref{TaylorDet}, we get
\begin{multline}\label{AppA4}
\frac{\mu}{2}\big(\text{tr} (\overline{\mathbf{1} + \mathbf{\Delta}}) -3 \big) =
\frac{\mu}{2} \big( (\det(\mathbf{1} + \mathbf{\Delta}) )^{-1/3} \text{tr} (\mathbf{1} + \mathbf{\Delta}) -3 \big) \\
=\frac{\mu}{2}\big( (1 - 1/3 \text{tr} \mathbf{\Delta} + 1/18 (\text{tr} \mathbf{\Delta})^2 +
1/6 \text{tr} (\mathbf{\Delta}^2) + O ({\Delta}^3) ) (3 + \text{tr} \mathbf{\Delta})  -3 \big) \\
=\frac{\mu}{4} \big( \text{tr} (\mathbf{\Delta}^2) - 1/3 (\text{tr} \mathbf{\Delta})^2 \big) + O ({\Delta}^3) =
\frac{\mu}{4} \text{tr} \big( (\mathbf{\Delta}^{\text{D}} )^2 \big) + O ({\Delta}^3).
\end{multline}
Finally, \eqref{EnerDist52} follows from \eqref{spec1}, using \eqref{AppA3} and \eqref{AppA4}.

\section*{Appendix B}

Let $ \mathbf{A}, \mathbf{B} \in M$ and
$\|\mathbf{A} - \mathbf{B}\| \rightarrow 0$.
Let us prove, for instance, that
\begin{equation}\label{AppB1}
\mathbf{B}^{-1} : (\mathbf{A} - \mathbf{B}) = O ( \|\mathbf{A} - \mathbf{B}\|^2 ).
\end{equation}
Indeed, since $\det (\cdot)$ is a smooth function, we have
\begin{equation}\label{AppB2}
\det (\mathbf{A}) = \det (\mathbf{B}) +
\frac{\displaystyle \partial \det (\mathbf{B})}{\displaystyle \partial \mathbf{B}} : (\mathbf{A} - \mathbf{B}) +
O ( \|\mathbf{A} - \mathbf{B}\|^2 ).
\end{equation}
Next, using the Jacobi formula, we get
\begin{equation}\label{AppB3}
\det (\mathbf{A}) = \det (\mathbf{B}) +
\det (\mathbf{B}) \mathbf{B}^{-\text{T}} : (\mathbf{A} - \mathbf{B}) +
O ( \|\mathbf{A} - \mathbf{B}\|^2 ).
\end{equation}
Finally, taking into account that $\det (\mathbf{A})= \det (\mathbf{B})=1$ and $\mathbf{B}^{-\text{T}} = \mathbf{B}^{-1}$,
we obtain \eqref{AppB1}.

\textbf{Remark.} Note that for the tangential space
$T_{\mathbf{B}} M$ to the manifold $M$ in $Sym$ we have
\begin{equation}\label{AppB4}
T_{\mathbf{B}} M = \{\mathbf{X} \in Sym \ | \ \mathbf{B}^{-1} : \mathbf{X} = 0 \}.
\end{equation}
Thus, relation \eqref{AppB1} implies that
 $\lim \limits_{ \mathbf{A} \rightarrow \mathbf{B} } \big( ( \mathbf{A} - \mathbf{B} )
 / \| \mathbf{A} - \mathbf{B} \| \big) \in T_{\mathbf{B}} M$ (if the limit exists).

\section*{Appendix C}

We need to construct a function $q(\theta)$ such that for small $\Delta$
\begin{equation}\label{AppC1}
q(\theta) \geq \Phi(\mathbf{C}^{(1)}_{\text{i}}, \mathbf{C}^{(2)}_{\text{i}} ) + O(\Delta) \quad
\text{for all} \ (\mathbf{C}^{(1)}_{\text{i}}, \mathbf{C}^{(2)}_{\text{i}} )
\in  S(\theta, \Delta).
\end{equation}

Let $(\mathbf{C}^{(1)}_{\text{i}}, \mathbf{C}^{(2)}_{\text{i}} )
\in  S(\theta, \Delta)$.
It follows from Appendix B, that
\begin{equation}\label{AppC2}
\big( \mathbf{C}^{(1)-1}_{\text{i}} - \mathbf{C}^{(2)-1}_{\text{i}} \big) : \mathbf{C}^{(1)}_{\text{i}} = O (\Delta^2).
\end{equation}
By $\mathbf{X}$ denote the orthogonal projection of $\mathbf{C}^{(1)-1}_{\text{i}} - \mathbf{C}^{(2)-1}_{\text{i}}$
on the tangential space $T_{\mathbf{C}^{(1)-1}_{\text{i}}} M$.
Using $\eqref{AppC2}$, we get for $\mathbf{X}$
\begin{equation}\label{AppC3}
\mathbf{X} = \mathbf{C}^{(1)-1}_{\text{i}} - \mathbf{C}^{(2)-1}_{\text{i}} -
\Big[ \big( \mathbf{C}^{(1)-1}_{\text{i}} - \mathbf{C}^{(2)-1}_{\text{i}} \big) : \mathbf{C}^{(1)}_{\text{i}}
\frac{1}{ \|\mathbf{C}^{(1)}_{\text{i}} \|^2 } \Big] \mathbf{C}^{(1)}_{\text{i}} =
\mathbf{C}^{(1)-1}_{\text{i}} - \mathbf{C}^{(2)-1}_{\text{i}} + O (\Delta^2).
\end{equation}
Moreover, since
$\text{tr} \big( \mathbf{X}  \mathbf{C}^{(1)}_{\text{i}}  \mathbf{X}  \big) = \| \mathbf{X} (\mathbf{C}^{(1)}_{\text{i}})^{1/2} \|^2$, we have
\begin{equation}\label{AppC32}
\frac{O (\Delta^3)}{\text{tr}
\big( \mathbf{X}  \mathbf{C}^{(1)}_{\text{i}}  \mathbf{X}  \big)} = O (\Delta).
\end{equation}
Substituting \eqref{AppC3} in \eqref{evaluat9} and taking \eqref{AppC32} into account, we obtain
\begin{equation}\label{AppC4}
\Phi(\mathbf{C}^{(1)}_{\text{i}}, \mathbf{C}^{(2)}_{\text{i}} ) =
\frac{-2 \big( \frac{(\mathbf{C}^{(1)-1}_{\text{i}})^{\text{D}}}{\|
(\mathbf{C}^{(1)-1}_{\text{i}})^{\text{D}} \|}
 : \mathbf{X} \big)
\big( \mathbf{1} : \mathbf{X} \big)}{\text{tr}
\big( \mathbf{X}  \mathbf{C}^{(1)}_{\text{i}}  \mathbf{X}  \big)} + O(\Delta).
\end{equation}

Thus, we define $q(\theta)$ as
\begin{equation}\label{AppC5}
q(\theta) := \max_{\| (\mathbf{C}^{(1)-1}_{\text{i}})^{\text{D}} \| \leq \theta}  \hat{q} (\mathbf{C}^{(1)}), \quad
\hat{q} (\mathbf{C}^{(1)}) := \max_{\mathbf{X} \in T_{\mathbf{C}^{(1)-1}_{\text{i}}} M}
\frac{-2 \big( \frac{(\mathbf{C}^{(1)-1}_{\text{i}})^{\text{D}}}{\| (\mathbf{C}^{(1)-1}_{\text{i}})^{\text{D}} \|}
 : \mathbf{X} \big)
\big( \mathbf{1} : \mathbf{X} \big)}{\text{tr}
\big( \mathbf{X}  \mathbf{C}^{(1)}_{\text{i}}  \mathbf{X}  \big)}.
\end{equation}

The function $\hat{q} (\mathbf{C}^{(1)}_{\text{i}})$ can be evaluated as follows.
First, for each $\mathbf{X}$ introduce
$\mathbf{Y} = \mathbf{X} ( \mathbf{C}^{(1)}_{\text{i}} )^{1/2}$.
Next, define a vector space $T := \{ \mathbf{Y} \in Sym \ | \
( \mathbf{C}^{(1)}_{\text{i}} )^{1/2} : \mathbf{Y} = 0 \}$. Thus,
\begin{equation}\label{AppC6}
\mathbf{X} \in T_{\mathbf{C}^{(1)-1}_{\text{i}}} M \quad \iff \mathbf{Y} \in T,
\end{equation}
\begin{equation}\label{AppC7}
\text{tr}
\big( \mathbf{X}  \mathbf{C}^{(1)}_{\text{i}}  \mathbf{X}  \big) = \| \mathbf{Y} \|, \quad
-2 \frac{(\mathbf{C}^{(1)-1}_{\text{i}})^{\text{D}}}{\| (\mathbf{C}^{(1)-1}_{\text{i}})^{\text{D}} \|}
 : \mathbf{X} = \mathbf{B}_1 : \mathbf{Y}, \quad
\mathbf{1} : \mathbf{X} = \mathbf{B}_2 : \mathbf{Y},
\end{equation}
where
\begin{equation}\label{AppC8}
\mathbf{B}_1 := -2 ( \mathbf{C}^{(1)}_{\text{i}} )^{-1/2}
\frac{(\mathbf{C}^{(1)-1}_{\text{i}})^{\text{D}}}{\|
(\mathbf{C}^{(1)-1}_{\text{i}})^{\text{D}} \|}, \quad
\mathbf{B}_2 := ( \mathbf{C}^{(1)}_{\text{i}} )^{-1/2}.
\end{equation}
Therefore,
\begin{equation}\label{AppC82}
\hat{q} (\mathbf{C}^{(1)}_{\text{i}}) = \max_{\mathbf{Y} \in T, \|\mathbf{Y}\| =1}
\Big[ \big( \mathbf{B}_1 : \mathbf{Y} \big)  \big( \mathbf{B}_2 : \mathbf{Y} \big) \Big].
\end{equation}
Next, we compute the orthogonal projections of $\mathbf{B}_1$ and $\mathbf{B}_2$ on $T$:
\begin{equation}\label{AppC9}
\mathbf{B}^0_{\text{k}} := \mathbf{B}_{\text{k}} -
\big( \mathbf{B}_{\text{k}} : ( \mathbf{C}^{(1)}_{\text{i}} )^{1/2}  \big)
( \mathbf{C}^{(1)}_{\text{i}} )^{1/2} \frac{1}{ \| ( \mathbf{C}^{(1)}_{\text{i}} )^{1/2} \|^2 }, \quad \text{k} \in \{1,2\}.
\end{equation}
Thus,
\begin{equation}\label{AppC10}
\hat{q} (\mathbf{C}^{(1)}_{\text{i}}) = \max_{\mathbf{Y} \in T, \|\mathbf{Y}\| =1}
\Big[ \mathbf{Y} :  \text{sym}\big( \mathbf{B}^0_1 \otimes \mathbf{B}^0_2 \big) : \mathbf{Y} \Big] =
\lambda_{\max} \big(\text{sym}\big( \mathbf{B}^0_1 \otimes \mathbf{B}^0_2 \big)\big),
\end{equation}
where $\lambda_{\max}\big(\text{sym}\big( \mathbf{B}^0_1 \otimes \mathbf{B}^0_2 \big)\big)$
is the maximal eigenvalue of the symmetric operator \\ $\text{sym}\big( \mathbf{B}^0_1 \otimes \mathbf{B}^0_2 \big): Sym \rightarrow Sym$.
Obviously, the same maximal eigenvalue has its restriction on $T^0 = \text{Span} \{\mathbf{B}^0_1, \mathbf{B}^0_2\}$.
It can be easily seen that
\begin{equation}\label{AppC102}
\text{sym}\big( \mathbf{B}^0_1 \otimes \mathbf{B}^0_2 \big) (\mathbf{B}^0_1) =
1/2 (\mathbf{B}^0_1 : \mathbf{B}^0_2) \mathbf{B}^0_1 + 1/2 (\mathbf{B}^0_1 : \mathbf{B}^0_1 ) \mathbf{B}^0_2,
\end{equation}
\begin{equation}\label{AppC103}
\text{sym}\big( \mathbf{B}^0_1 \otimes \mathbf{B}^0_2 \big) (\mathbf{B}^0_2) =
1/2 (\mathbf{B}^0_2 : \mathbf{B}^0_2) \mathbf{B}^0_1 + 1/2 (\mathbf{B}^0_1 : \mathbf{B}^0_2 ) \mathbf{B}^0_2.
\end{equation}
Therefore, the matrix of the restricted operator with
respect to the basis $ \{\mathbf{B}^0_1, \mathbf{B}^0_2\}$
has the following form
\begin{equation}\label{AppC11}
A := \frac{1}{2} \left(
\begin{array}{cc}
\mathbf{B}^0_1 : \mathbf{B}^0_2 \ & \ \mathbf{B}^0_2 : \mathbf{B}^0_2   \\
\mathbf{B}^0_1 : \mathbf{B}^0_1 \ & \ \mathbf{B}^0_1 : \mathbf{B}^0_2
\end{array}
\right).
\end{equation}

Both eigenvalues of $A$ are real, since $A$ represents a symmetric tensor.
Finally,
\begin{equation}\label{AppC12}
\hat{q} (\mathbf{C}^{(1)}_{\text{i}}) = \lambda_{\max}\big(\text{sym}\big( \mathbf{B}^0_1 \otimes \mathbf{B}^0_2 \big)\big) = \lambda_{\max}(A).
\end{equation}

\begin{figure}\centering
\psfrag{A}[m][][1][0]{$\theta$}
\psfrag{B}[m][][1][0]{$q$}
\scalebox{0.8}{\includegraphics{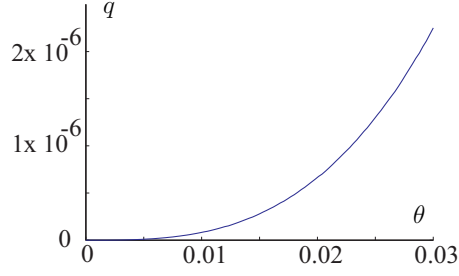}}
\caption{Function $q(\theta)$. \label{fig5}}
\end{figure}

Note that $\hat{q} (\mathbf{C}^{(1)}_{\text{i}})$
is a continuous function of $\mathbf{C}^{(1)}_{\text{i}}$.
Therefore, the maximum
$q(\theta) = \max \limits_{\| (\mathbf{C}^{(1)-1}_{\text{i}})^{\text{D}} \| \leq \theta}  \hat{q} (\mathbf{C}^{(1)})$
is well defined. We compute it by the brutal force method. Moreover, the following
parametrization can be used to simplify the computations.
For any tensor $\mathbf{C}^{(1)}_{\text{i}}$ there exists a cartesian coordinate system
and real numbers $\lambda^1, \lambda^2 > 0$ such that
the matrix of $\mathbf{C}^{(1)}_{\text{i}}$ takes the diagonal form
$\text{diag}(\lambda^1, \lambda^2, 1/( \lambda^1 \lambda^1) )$.
The function $q(\theta)$ is plotted on the figure \ref{fig5} for $0 \leq \theta \leq 0.03$.

\section*{Acknowledgements}

This research was supported by German National Science Foundation (DFG) within the
collaborative research center SFB 692 "High-strength aluminium based light weight materials for
reliable components".




\end{document}